\documentclass[12pt,reqno]{amsart}
\usepackage{amsthm,amsfonts,amssymb,amscd,amsmath,euscript}
\usepackage[all]{xy}
\newtheorem{theorem}[subsection]{Theorem}

\newtheorem{proposition}[subsection]{Proposition}

\newtheorem{lemma}[subsection]{Lemma}
\newtheorem{corollary}[subsection]{Corollary}

\theoremstyle{definition}
\newtheorem{definition}[subsection]{Definition}

\newtheorem{proposition-definition}[subsection]{Proposition-Definition}

\theoremstyle{remark}
\newtheorem{remark}[subsection]{Remark}

\newcommand{\connsum}[1]{{\raisebox{-2.ex}{$\stackrel{{\textstyle \#}}
{{\scriptstyle (#1)}}$}}}

\newcommand{\ccd}{{\begin{picture}(2,2)
\put(1,1){\circle*{0.7}}\end{picture}}}

\newcommand{\smallcup}{{\raisebox{0.2 ex}{$\scriptscriptstyle \cup$}}}

\newcommand{\dual}{{\scriptscriptstyle \vee}}
\newcommand{\sdual}{{\,\check{}}}

\newcommand{\Aut}{\operatorname{Aut}\nolimits}

\newcommand{\Cl}{\operatorname{Cl}}
\newcommand{\codim}{\operatorname{codim}\nolimits}

\newcommand{\Ext}{\operatorname{Ext}\nolimits}
\newcommand{\EXT}{{{\mathcal E}xt\:}}
\newcommand{\Fix}{\operatorname{Fix}}
\newcommand{\Fl}{\operatorname{Fl}}

\newcommand\Hilb{{\operatorname{Hilb}\nolimits}}

\newcommand{\HOM}{{{\mathcal H}om\:}}
\newcommand{\Hom}{\operatorname{Hom}\nolimits}
\newcommand{\id}{\operatorname{id}\nolimits}
\newcommand{\im}{\operatorname{im}\nolimits}
\newcommand{\Indet}{\operatorname{Indet}\nolimits}
\newcommand{\length}{\operatorname{length}\nolimits}

\newcommand{\norm}{{\operatorname{norm}\nolimits}}
\newcommand\ns{{\operatorname{ns}\nolimits}}

\newcommand{\PAut}{\operatorname{PAut}\nolimits}

\newcommand{\Pic}{\operatorname{Pic}\nolimits}
\newcommand{\PPic}{\operatorname{\overline{Pic}\:}\nolimits}
\newcommand{\Proj}{\operatorname{Proj}}
\newcommand{\Prym}{\operatorname{Prym}\nolimits}
\newcommand{\PPrym}{\operatorname{\overline{Prym}\:}\nolimits}
\newcommand{\pr}{\operatorname{pr}\nolimits}

\newcommand{\rk}{\operatorname{rk}\nolimits}
\newcommand{\Tr}{\operatorname{Tr}}
\newcommand{\Sing}{\operatorname{Sing}}

\newcommand{\Spec}{\operatorname{Spec}}

\newcommand{\Supp}{\operatorname{Supp}}

\newcommand{\CC}{{\mathbb C}}

\newcommand{\PP}{{\mathbb P}}
\newcommand{\QQ}{{\mathbb Q}}
\newcommand{\RR}{{\mathbb R}}

\newcommand{\ZZ}{{\mathbb Z}}

\newcommand{\GGB}{{\mathbf G}}

\newcommand{\OOO}{{\mathcal O}}

\newcommand{\JJJ}{{\mathcal J}}

\newcommand{\EEE}{{\mathcal E}}

\newcommand{\LLL}{{\mathcal L}}
\newcommand{\FFF}{{\mathcal F}}
\newcommand{\CCC}{{\mathcal C}}

\newcommand{\NNN}{{\mathcal N}}
\newcommand{\MMM}{{\mathcal M}}
\newcommand{\QQQ}{{\mathcal Q}}

\newcommand{\GGGG}{\boldsymbol{\mathcal G}}
\newcommand{\LLLL}{\boldsymbol{\mathcal L}}
\newcommand{\PPPP}{\boldsymbol{\mathcal P}}

\renewcommand{\bar}[1]{\overline{#1}}

\newcommand\alp{\alpha}

\newcommand\G{\Gamma}

\newcommand\eps{\epsilon}

\renewcommand\phi{\varphi}

\newcommand\ups{\upsilon}

\newcommand{\equi}{\: \Longleftrightarrow\: }
\newcommand{\isoto}{{\lra\hspace{-1.3 em}
\raisebox{ 0.6 ex}{$\textstyle\sim$}\hspace{0.8 em}}}
\newlength{\rrrr}
\newcommand{\into}{\hookrightarrow}
\newcommand{\onto}{\twoheadrightarrow}

\newcommand\lra{{\longrightarrow}}

\newcommand\rar{\rightarrow}
\newcommand\lrdash{\:
\xymatrix@1{\ar@{-->}[r]&}\:
}
\newcommand{\lrdashar}[1]{\:
\xymatrix@1{\ar@{-->}[r]^{#1}&}\:
}
\newcommand{\sendsto}[1]{\:
\xymatrix@1{\ar@{|->}[r]^{#1}&}\:
}

\newcommand\empt{\varnothing}

\author{D. Markushevich}
\address{\scriptsize D. Markushevich:  
Math\'ematiques - b\^{a}t.M2, Universit\'e Lille 1, 
F-59655 Villeneuve d'Ascq Cedex, France}
\email{markushe@math.univ-lille1.fr}

\author{A.S.~Tikhomirov}

\address{\scriptsize   A.S.~Tikhomirov:
  Department of Mathematics\\
  State Pedagogical University\\
  Respublikanskaya Str. 108
\newline 150 000 Yaroslavl, Russia}
\email{tikhomir@yaroslavl.ru}


\title{New symplectic $V$-manifolds of dimension four
via the relative compactified Prymian}

\begin{document}

\begin{abstract}
Three 
new examples of 4-dimensional irreducible symplectic
$V$-manifolds are constructed.
Two of them are relative compactified
Prymians of a family of genus-3 curves with involution,
and the third one is obtained from a Prymian by Mukai's flop.
They have the same singularities as two of Fujiki's examples,
namely, 28 isolated singular points analytically equivalent
to the Veronese cone of degree 8, but a different Euler number.
The family of curves used in this construction forms
a linear system on a K3 surface with involution.
The structure morphism of both Prymians to the base of the family
is a Lagrangian fibration in abelian surfaces
with polarization of type (1,2). No example of such fibration
is known on {\em nonsingular} irreducible symplectic varieties.
%
\end{abstract}

\maketitle



\setcounter{section}{-1}

\section{Introduction}

Historically, the first constructions of nontrivial
compact K\"ahler holomorphically symplectic
varieties $Y$ of dimension $>2$ belong to
Beauville \cite{Beau-1} and Fujiki \cite{F}. Fujiki's notion of
nontriviality means
that $Y$ is not obtained as a finite quotient
from a product of a complex torus with
symplectic surfaces. Fujiki constructed one nonsingular example in dimension 4,
the blowup $S^{[2]}$ of the diagonal in the symmetric
square $S^{(2)}$ of a K3 surface $S$, and his other examples
are 4-dimensional $V$-manifolds, that is varieties having
finite quotient singularities.

Beauville
\cite{Beau-1}, \cite{Beau-2} constructed two deformation classes of
nonsingular irreducible symplectic manifolds
in all even dimensions $2n$. 
Here a manifold is called {\em irreducible symplectic}
if it is simply connected and has a unique
symplectic structure up to proportionality; this is equivalent to
Fujiki's condition of nontriviality at least in the category of
nonsingular symplectic varieties.
The Beauville's examples are: 1)~$S^{[n]}=\Hilb^n(S)$, the Hilbert scheme of
$0$-dimensional subschemes of length $n$ in a K3 surface $S$, and  2) $K_n(A)$,
the generalized Kummer variety associated to an abelian surface $A$.
The latter is defined as the fiber of the summation map $A^{[n+1]}\rar A$.

Mukai \cite{Mu-1} showed that the moduli spaces of semistable sheaves
on a K3 or abelian surface are symplectic. According to
\cite{Hu-1}, \cite{Hu-2}, \cite{O'G-1} and \cite{Y}, whenever
such a moduli space is nonsigular, it is deformation equivalent
to $S^{[n]}$ or $K_n(A)\times T$ with $T=A$ or $A\times A$. Thus, for years,
two Beauville's examples provided the only known moduli components
of irreducible symplectic manifolds, until O'Grady \cite{O'G-2},
\cite{O'G-3} constructed two essentially new such manifolds.
They are obtained as symplectic desingularizations of singular moduli
spaces of semistable sheaves. The first one is associated to a K3
surface and is of dimension 10, the second one is associated to
an abelian surface and is of dimension 6. It is still unknown
whether there exist irreducible symplectic 4-folds that are
not deformation equivalent to $S^{[2]}$ or  $K_2(A)$. O'Grady
studies in \cite{O'G-4}, \cite{O'G-5} the irreducible symplectic 4-folds
whose intersection 4-linear form on $H^2$ is isomorphic to that of
$S^{[2]}$ and conjectures that they are deformation equivalent to
$S^{[2]}$.

The results of \cite{KL}, \cite{KLS}, \cite{LS}, \cite{CK-1},
\cite{CK-2} show that, informally speaking, no new examples
of {\em nonsingular} irreducible symplectic manifolds can be
obtained by the method of \cite{O'G-2},
\cite{O'G-3}. More precisely, for any singular moduli
space $\MMM$ of semistable sheaves on a K3 surface, either $\MMM$
has no symplectic resolution, or such a resolution exists and
up to deformations coincides with one of the known examples:
Beauvilles's or O'Grady's. A weaker result, concerning only
rank-2 sheaves, is obtained for moduli of sheaves on abelian surfaces.
 
Thus the problem of extending the very short list of known
deformation classes of irreducible symplectic {\em manifolds}
is very hard. Leaving aside this hard problem,
we turn back to the original setting of Fujiki, who considered
symplectic {\em $V$-manifolds}. All of his examples, up to
deformation of a complex structure, are
partial resolutions of finite quotients
of the products of two symplectic surfaces. 


In the present article, we provide a new construction 
of irreducible symplectic $V$-manifold of dimension 4,
the relative compactified Prym variety
of some family of curves with involution. The fibration in Prym surfaces
is Lagrangian. 

Many features of the theory of irreducible
symplectic manifolds are very similar to those of K3 surfaces,
and in view of this similarity, the manifolds with a Lagrangian fibration
constitute an important class of irreducible
symplectic manifolds which is an analog of the class
of K3 surfaces with an elliptic pencil.
Earlier examples of Lagrangian fibrations on irreducible
symplectic manifolds were constructed in \cite{Beau-3}, \cite{D},
\cite{HTsch-1}, \cite{HTsch-2}, \cite{IR}, \cite{S-2}. 
There are more examples
if we relax the hypothesis that the fibration map is a regular morphism,
but admit {\em rational} Lagrangian fibrations \cite{Mar-2}, \cite{Gu}.

By Liouville's Theorem, the general fiber of these Lagrangian fibrations
is an abelian variety. It turns out, that 
in all of the known examples it is 
either the Jacobian of a curve,   or a $(1,\ldots,1,k)$-polarized
abelian variety with $k\geq 3$.
The first possibility occurs
for Lagrangian fibrations $f:Y\rar B$ with $Y$ deformation
equivalent to $S^{[n]}$, and the second one
for $Y$ birational to $K_n(A)$, where $A$ is an abelian
surface of polarization $(1,k)$ (see \cite{S-1}, Remark 3.9). 

Thus there
are no examples of Lagrangian fibrations on
irreducible symplectic 4-folds with $(1,2)$-polarized
abelian surfaces as fibers. On the other hand,
there are classical integrable systems integrated on Prym surfaces
of such polarization, for example, the complexified Kowalevski top
\cite{HvM}. However, the corresponding symplectic manifolds,
which are the (complexified) phase spaces of these systems,
are always rational, or at least unirational, and hence they are
very far from having a symplectic compactification, neither nonsingular,
nor in the class of $V$-manifolds. The phase space of
the Kowalevski top is identified with
the relative Prym variety $\Prym^k(\CCC,\tau)$ of a family $\CCC/\PP^2$ 
of genus-3 curves endowed with an involution $\tau$ such that
the quotients by $\tau$ form a family of elliptic curves.

In the present paper, we use this idea in taking for $\CCC/\PP^2$
the family of $\tau$-invariant members of 
a linear system $|H|$ of genus-3 curves on a K3 surface $S$
with an involution $\tau$.
In order that the construction might work, $\tau$ should leave
the symplectic form $\omega\in H^0(S,\Omega^2_S)$ anti-invariant.

We denote by $\PPrym^{k,\kappa}(\CCC,\tau)$ ($k\in\ZZ$) the relative compactified 
Prym variety defined as a connected component of the fixed locus of some involution $\kappa$
in the relative compactified  Picard
variety $\PPic^k(|H|)$. The latter is a compactification
of the relative Picard variety $\Pic^k(|H|)$ parametrizing divisor classes of degree $k$ on the
curves from the complete linear system $|H|$ of an ample genus-3 curve $H$.
The compactification depends on the choice of a polarization on $S$, which
we fix once and forever to be $H$.
There are at most 4 non-isomorphic compactified Picard
varieties, corresponding to $k=0,1,2,3$, and $\PPic^k(|H|)$ is birational
to $\PPic^{k+2}(|H|)$, so there are at most two nonbirational ones. The definiton
of $\kappa$ depends on some arbitrary choices if $k=1$ or $3$, but is
canonical for even $k$, so for even $k$, we suppress the superscript
$\kappa$ from the notation. We work out in full detail the case of even $k=2m$,
proving that $\PPPP^{2m}=\PPrym^{2m}(\CCC,\tau)$ is an irreducible symplectic
$V$-manifold. As $\PPPP^{2m}\simeq \PPPP^{2m+4}$,
there are at most two non-isomorphic compactified Prymians,
$\PPPP^{0}$ and $\PPPP^{2}$. We do not know whether they
are really non-isomorphic. 

The Prymian $\PPPP^{0}$ of degree 0 contains a family of groups,
hence has the zero section whose image $\Pi$ is isomorphic to the
base of the family, that is to $\PP^2$. Performing Mukai's flop
with center $\Pi$, we obtain a third 4-dimensional symplectic $V$-manifold
$M'$ (see Corollary \ref{M-prime}). 
It has the same topological invariants and Hodge numbers as $\PPPP^0$,
but is conjectured to be non-isomorphic neither to $\PPPP^{0}$,
nor to $\PPPP^{2}$. We also identify $M'$
as a partial desingularization of
the quotient of $S^{[2]}$ by a symplectic involution.

In the case of odd degree $k$, it is unlikely that $\PPrym^{k,\kappa}(\CCC,\tau)$
is symplectic. We show that it is not symplectic for one of the possible choices of
$\kappa$, for which it contains a $3$-dimensional rational variety, see
Remark~\ref{odd-degree}.\hyphenpenalty=10000

\bigskip

{\sc Acknowledgements.} The authors thank Eyal Markman and Christoph
Sorger for discussions. They acknowledge with pleasure the hospitality of
the Max-Planck-Institut f\"ur Mathematik in Bonn, where was done
a part of the work on the present paper. The second author also
thanks for hospitality the University of Lille I.

\section{Definition and basic properties of varieties $\PPic^k(|H|)$}
\label{section1}

Let $X$ be a Del Pezzo surface of degree 2 obtained as a
double cover of $\PP^2$ branched in a generic quartic
curve $B_0$, $\mu:X\lra\PP^2$ the double cover map,
$B=\mu^{-1}(B_0)$ the ramification curve. Let $\Delta_0$
be a generic curve from the linear system $|-2K_X|$,
$\rho:S\lra X$ the double cover branched in $\Delta_0$
and $\Delta=\rho^{-1}(\Delta_0)$.
Then $S$ is a K3 surface, and $H=\rho^*(-K_X)$ is a degree-4 ample
divisor class on $S$, which we fix once and forever
as a polarization of $S$. We will denote by $\iota$
(resp. $\tau$) the Galois involution of the double cover $\mu$
(resp. $\rho$).

The plane quartic $B_0$
has 28 bitangent lines $m_1,\ldots ,m_{28}$,
and $\mu^{-1}(m_i)$ is the union
of two rational curves $\ell_i\cup\ell'_i$ meeting in 2 points.
The 56 curves $\ell_i,\ell'_i$ are all the {\em lines}
on $X$, that is, curves of degree 1 with respect to $-K_X$.
Further, the curves
$C_i=\rho^{-1}(\ell_i)$, $C'_i=\rho^{-1}(\ell'_i)$ are
{\em conics} on $S$, that is, curves of degree 2
with respect to $H$. Each pair $C_i,C'_i$ meets in
4 points, thus forming a reducible curve of arithmetic
genus 3 belonging to the linear system $|H|$.
Throughout the paper we assume that $B_0$, $\Delta_0$
are sufficiently generic. This implies, in particular, that
each line $\ell_i$ meets only one of the two lines $\ell_j,\ell'_j$
for $j\neq i$. 

\begin{lemma}\label{pairs-of-conics}
The linear system $|H|$ is very ample and embeds
$S$ into $\PP^3$ as a quartic surface. Every curve
in $|H|$ is reduced, and the only reducible members
of $|H|$ are the $28$ curves $\G_i=C_i+C'_i$, $i=1,\ldots,28$.
\end{lemma}

\begin{proof}
A generic curve in $|H|$ is isomorphic to a plane quartic, for
$\Delta\in |H|$ and $\Delta\simeq\Delta_0\simeq\mu(\Delta_0)\in
|\OOO_{\PP^2}(4)|$. According to Saint-Donat's description
of ample linear systems on K3 surfaces \cite{SD} (see also \cite{Mor},
Section 6), $|H|$ is very ample and embeds $S$ into $\PP^3$.

Let $C\in |H|$ be reducible or non-reduced. Then the
same is true for $\underline{C}=\rho_* (C)\in |-2K_X|$.
If $\underline{C}=D_1+D_2$ is reducible and reduced,
then we have the following 3 possibilities: (a) $\mu(D_1)$
is a line and $\mu (D_2)$ is a cubic; (b) $\mu (D_1),\mu (D_2)$
are conics; (c) $\mu$ is of degree 2 over one or both
components $D_i$.

In the cases (a) and (b), $\mu^{-1}\mu (D_i)$ decomposes
into two components $D_i$ and $\iota (D_i)$ for both $i=1,2$.
Hence $\mu (D_i)$ are totally tangent to $B_0$. Hence the family of
such curves $D_1+D_2$ is 3-dimensional in case (a)
and 2-dimensional in case (b). Similarly, $\rho ^{-1}(D_i)$
is the union of two components permuted by $\tau$,
and $D_i$ are totally tangent to $\Delta_0$ for both $i=1,2$.
This is impossible for generic $B_0$, $\Delta_0$
by dimension reasons.

In the case (c), let $E=\rho^{-1}(D_i)$, where $D_i$
is the component mapped onto the line $\mu (D_i)$
with degree 2. Then $E^2=0$, $H\cdot E=2$, which is impossible
by loc. cit., for then every smooth member of $|H|$ should be hyperelliptic.

By a similar argument, one can eliminate the case when $\rho (C)$
is reducible and $\rho_*(C)$ is non-reduced. Thus, the only
remaining case is when $\deg \rho|_C=2$ and $\deg \mu|_{\rho (C)}=2$,
in which $\mu\rho(C)$ is a bitangent to~$B_0$.
\end{proof}

Mukai \cite{Mu-2} has endowed the integer cohomology
$H^*(Y)$ of a K3 surface $Y$ with the following bilinear form:
$$
\langle(v_0,v_1,v_2),(w_0,w_1,w_2)\rangle=v_1\smallcup w_1-
v_0\smallcup w_2-v_2\smallcup w_0,
$$
where $v_i,w_i\in H^{2i}(X)$. For a sheaf $\FFF$ on $Y$,
the {\em Mukai vector} is $v(\FFF)=(\rk \FFF, c_1(\FFF),
\chi (\FFF)-\rk \FFF)\in H^*(Y)$, where $H^4(Y)$ is naturally identified
with $\ZZ$. 
We refer to \cite{Sim} or to \cite{HL} for the definition
and the basic properties of the Simpson (semi-)stable sheaves.
Let $M_Y^{H,s}(v)$ (resp. $M_Y^{H,ss}(v)$) denote
the moduli space of Simpson stable (resp. semistable) \cite{Sim}
sheaves $\FFF$
on $Y$ with respect to an ample class $H$ with
Mukai vector $v(\FFF)=v$.
According to Mukai (\cite{Mu-1}, \cite{Mu-2}, see also \cite{HL}),
$M_Y^{H,s}(v)$, if nonempty,
is smooth of dimension $\langle v,v\rangle +2$ and carries
a holomorphic symplectic structure.
We will study the moduli space $\MMM^k=M^{H,ss}_S(v)$
on the above special K3 surface $S$
with Mukai vector $v=(0,H,k-2)$. 

\begin{proposition}\label{basic}
(i) $\MMM=\MMM^k$ is an irreducible
projective variety of dimension~$6$. 
The open part $\MMM^{*}=M^{H,s}_S(0,H,k-2)$ corresponding to the stable
sheaves is contained in the smooth locus of $\MMM$ and is a holomorphically
symplectic manifold with symplectic form $\alp\in H^0(\MMM^*,\Omega^2)$
induced by the Yoneda pairing
$$
\alp_{[\LLL]}:\Ext^1(\LLL,\LLL)\times
\Ext^1(\LLL,\LLL)\lra \Ext^2(\LLL,\LLL)\xrightarrow{\Tr}H^2(S,\OOO_S)\simeq\CC,
$$
where $[\LLL]\in \MMM^*$ is the class of a stable sheaf $\LLL$
and the tangent space $T_{[\LLL]}\MMM^*$ is identified
with $\Ext^1(\LLL,\LLL)$.

(ii) $\MMM^k$ parametrizes the $S$-equivalence classes
of pure $1$-dimensional sheaves $\LLL$ whose supports
are curves from the linear system $|H|$ and
such that $\LLL|_C$ is a torsion free  $\OOO_C$-module of rank $1$
with $\chi (\LLL)=k-2$, where $C=\Supp\LLL$. In the case when $\LLL$ is
invertible as a sheaf on its support, the condition
$\chi (\LLL)=k-2$ is equivalent to saying that $\deg \LLL=k$.

(iii) For any $k\in\ZZ$, $\MMM^k\simeq\MMM^{k+4}$. For odd $k$, any semistable sheaf
from $\MMM^k$ is stable, so $\MMM^k$ is nonsingular.
For even $k$, $\MMM^k$ contains exactly
$28$ $S$-equivalence
classes of strictly semistable sheaves.
Each of them is the class of the sheaf
$\OOO_{C_i}(\frac{k-4}{2}pt)\oplus\OOO_{C'_i}(\frac{k-4}{2}pt)$
($i=1,\ldots,28$), where $pt$ stands for the class of a point
on either one of the conics $C_i,C'_i$.

\end{proposition}

\begin{proof}
(i) The projectivity of $\MMM^k$ follows by Theorem 1.21 of \cite{Sim}.
The stable sheaves being simple, the remaining assertions
follow by Theorem~0.1 of \cite{Mu-1}. 

(ii) If $[\LLL]\in\MMM^k$, then $v(\LLL)=(0,H,k-2)$, so $\LLL$
is, by definition, an equidimensional torsion sheaf
with $c_1(\LLL)=H$ and $\chi(\LLL)=k-2$.
Hence the support of $\LLL$ is a curve from $|H|$ and the rank
of $\LLL$ is 1 at the generic points of all the components of $H$.
It is torsion-free when considered as a sheaf on $C$ since
it is equidimensional.

(iii) The isomorphism $\MMM^k\rar\MMM^{k+4}$ is given by $[\LLL]\mapsto [\LLL(1)]$
for all $[\LLL]\in\MMM^k$. Further, if
$[\LLL]\in\MMM^k$ and $C=\Supp\LLL$ is an integral curve, then
any rank-1 torsion-free sheaf on $C$ is stable according to Simpson's
definition, whether it is considered as a sheaf on $C$ or on $S$,
for it has no proper 1-dimensional subsheaves. By (i), this implies
that $[\LLL]$ is a smooth point of $\MMM^k$. Now suppose that
$C$ is not integral. By Lemma \ref{pairs-of-conics}, $C$
is one of the curves $C_i+C'_i$. Hence the only possibility
for a strictly semistable sheaf is to be the central term of
an extension
\begin{equation}\label{triple}
0\lra\FFF\lra\LLL\lra\FFF'\lra 0,
\end{equation}
where $\FFF,\FFF'$ are pure 1-dimensional,
$\Supp\FFF=C_i$, $\Supp\FFF'=C'_i$, or vice versa,
and $\chi (\FFF (n))=\chi (\FFF' (n))$ for $n\gg 0$.
Hence $\FFF,\FFF'$ are invertible on their supports,
and $\chi (\FFF)=\chi (\FFF')=\frac{k-2}{2}$, which implies
that $k$ is even and $\FFF\simeq\OOO_{C_i}(\frac{k-4}{2}pt)$, 
$\FFF'\simeq\OOO_{C'_i}(\frac{k-4}{2}pt)$.
\end{proof}

\begin{definition}\label{ISVM}
(i) A $V$-manifold is an algebraic variety having at worst finite quotient
singularities. We reserve the term ``manifold'' for nonsingular
algebraic varieties.

(ii) A symplectic variety is a normal algebraic variety $Y$ such that its
nonsingular
locus $Y_{\ns}$ has a symplecitc structure,
that is a $2$-form $\omega\in H^0(Y_{\ns},\Omega^2_{Y_{\ns}})$
 which is closed and everywhere
nondegenerate on~$Y_{\ns}$. The nondegeneracy means that $\omega^{\wedge\frac{1}{2}\dim Y}$
has no zeros on $Y_{\ns}$. If $Y$ is nonsingular, we also call it
a symplectic manifold.

(iii) A closed irreducible subvariety $W\subset Y$ of a symplectic variety $Y$
endowed with a symplectic structure $\omega$
is called Lagrangian if $\dim W=\frac{1}{2}\dim Y$, 
$W_0:=Y_{\ns}\cap W_{\ns}\neq\empt$
and $\omega|_{W_0}\equiv 0$.

(iv) A symplectic manifold (or $V$-manifold) $Y$ is said to be
irreducible symplectic if $Y$ is complete, simply connected, and
$h^0(Y,\Omega^2_Y)=1$.


(v) A morphism $f:Y\rar B$ from a symplectic variety
of dimension $2n$ to another variety $B$ of dimension $n$
is called a Lagrangian fibration if it is surjective and if
the generic fiber of $f$ is a connected Lagrangian subvariety of $Y$.
\end{definition}

\begin{proposition}\label{fibers-of-f}
In the above notation, the map
$f:\MMM^k\rar |H|\simeq\PP^3$ sending $[\LLL]\in\MMM^k$
to the curve $C_\LLL=\Supp\LLL\in |H|$ is a Lagrangian fibration.
The following properties are verified:

(i) If $C\in |H|$ is smooth, then the fiber $f^{-1}(\{C\})$ is canonically
isomorphic to $\Pic^k(C)$. Here $\{C\}$ denotes the point of the projective
$3$-space $|H|$ representing the curve $C$.
Further, if $U\subset|H| $ is the open set parametrizing
integral curves, $U=|H|\setminus \{\{\G_1\},\ldots,\{\G_{28}\}\}$, then the
restriction $f_U:f^{-1}(U)\rar U$ of $f$ over $U$ 
is identified with the relative compactified  Picard variety of 
Altman--Kleiman.

(ii) Let $\Gamma=C+C'$ be one of the reducible curves
$\G_i=$ ($i=1,\ldots,28$).

If $k$ is even,  then  $f^{-1}(\{C\})$ is the
union of three $3$-dimensional rational components
$\bar{J}^{\frac{k-2}{2},\frac{k+2}{2}},\bar{J}^{\frac{k}{2},\frac{k}{2}},\bar{J}^{\frac{k+2}{2},\frac{k-2}{2}}$. 

If $k$ is odd, then $f^{-1}(\{C\})$ is the
union of four $3$-dimensional rational components
$\bar{J}^{\frac{k-3}{2},\frac{k+3}{2}},\ldots,\bar{J}^{\frac{k+3}{2},\frac{k-3}{2}}$.

Each $\bar{J}^{d,d'}=\bar{J}^{d,d'}(\G)$ contains
an open subset $J^{d,d'}=J^{d,d'}(\G)\simeq(\CC^*)^3$ parametrizing
the invertible $\OOO_{\G}$-modules $\LLL$ such that
$\deg\LLL|_{C}=d$, $\deg\LLL|_{C'}=d'$.

(iii) Let $q$ be one of the $56$ conics $C_i,C'_i$  ($i=1,\ldots,28$).
Then $\MMM^k$ is birational to $\MMM^{k+2}$ via the map
$\psi: [\LLL]\mapsto [\LLL (q)]$. Let us set $q=C_i$, and fix
the notation for the $56$ conics in such a way that
$C_i\cap C'_j=\empt$ and $C_i\cap C_j= 2$ points
for all $j\neq i$.

If $k$ is even,  then the indeterminacy locus of $\psi$
is given by the formula
$$\Indet(\psi)=f^{-1}(\{\G_i\})\cup
\bigcup_{j\neq i}\bar{J}^{\frac{k+2}{2},\frac{k-2}{2}}(\G_j).$$ 

If $k$ is odd,  then the indeterminacy locus of $\psi$ is
\begin{multline}\label{indetpsi}
\Indet(\psi)=\bar{J}^{\frac{k-3}{2},\frac{k+3}{2}}(\G_i)\cup 
\bar{J}^{\frac{k-1}{2},\frac{k+1}{2}}(\G_i)\cup\\ 
\bar{J}^{\frac{k+1}{2},\frac{k-1}{2}}(\G_i)
\cup\bigcup_{j\neq i}\bar{J}^{\frac{k+3}{2},\frac{k-3}{2}}(\G_j).
\end{multline}

The formulas for $\Indet(\psi)$ in the case when $q=C'_i$ are obtained by
replacing all the $\bar{J}^{m,n}$ by $\bar{J}^{n,m}$.
\end{proposition}

\begin{proof}
(i) The map $f:\MMM^k\lra |H|$ can be defined as a map from 
the moduli functor of sheaves
on $S$ to the Hilbert functor of curves on $S$, using
the $0$-th Fitting ideal of a torsion sheaf, and it
obviously commutes with base change and descends to the
schemes $\MMM^k, \Hilb_S$ representing these functors.

Let $U\subset|H|$ be the complement of the 28 reducible
curves, and $\phi:\CCC_U\lra U$ the universal curve of the linear system
$|H|$, restricted over $U$. Every fiber $C_t=\phi^{-1}(t)$
for $t\in U$ is an integral curve, so the Altman-Kleiman
relative compactified Jacobian $\bar{J}^k\phi:\bar{J}^k(\CCC_U/U)\lra U$ is defined
\cite{AK}, which is the relative moduli space parametrizing
the isomorphism classes of
degree-$k$ torsion-free rank-$1$ sheaves on the fibers of
$\phi$ ($d\in\ZZ$). It commutes with base change, so
$(\bar{J}^k\phi)^{-1}(t)=\bar{J}^k(C_t)$. Since the curves
$C_t$ lie on a smooth surface, they may have only planar
singularities. Then by \cite{AIK}, $\bar{J}^k(\CCC_U/U), \bar{J}^k(C_t)$
are reduced and irreducible and are compactifications
of the Picard schemes $\Pic^k(\CCC_U/U)$, resp. $\Pic^k(C_t)$.
By the universal property of
moduli spaces, there is a natural morphism $\bar{J}^k(\CCC_U/U)\lra\MMM^k$
which is bijective onto its image, equal to $f^{-1}(U)$. As
$f^{-1}(U)$ is nonsingular, $\bar{J}^k(\CCC_U/U)$ is nonsingular
($=$ smooth over $\CC$) as well and the last map is an isomorphism
identifying $f_U=f|_{f^{-1}(U)}$ with $\bar{J}^k\phi$. By \cite{Beau-3},
$f$ is Lagrangian for $k=3$, the genus of the curves from
$|H|$. As $\Pic^k(\CCC_U/U)$ is a torsor under $\Pic^0(\CCC_U/U)$ 
in the \'etale topology, $f$ is
Lagrangian for any $k$ by \cite{MarT}, Lemma 5.7.

(ii) Let $k$ be even; the case of odd $k$ is completely similar. 
We are to show that the special fibers
$f^{-1}(t_i)$, where
$\{t_1,\ldots, t_{28}\}= |H|\smallsetminus U$, 
are unions of 3 components. Let $\G_i$ be represented by $t_i$, 
and look again
at the exact triple (\ref{triple}), but now
$\FFF,\FFF'$ are invertible on their supports
with $\chi(\FFF)\leq 0\leq \chi(\FFF')$. Let $C_i\cap C'_i=
\{z_1,\ldots,z_4\}$. In each point $z_k$, the stalk of the sheaf
$\EXT^1_{\OOO_S}(\FFF,\FFF')$ is a 1-dimensional vector space
$\CC_{z_k}$, so, locally at $z_k$, there are only two non-isomorphic
extensions: $\LLL_{z_k}\simeq \OOO_{\G_i,z_k}$ (the non-trivial extension)
and $\LLL_{z_k}\simeq \OOO_{C_i,z_k}\oplus\OOO_{C'_i,z_k}$ (the trivial
one). We have 
\begin{equation}\label{ext-one}
\Ext^1(\FFF,\FFF')= H^0(\EXT^1_{\OOO_S}(\FFF,\FFF'))
\simeq \bigoplus_{k=1}^4\CC_{z_k},
\end{equation}
so that every $\xi\in \Ext^1(\FFF,\FFF')$ can be viewed
as a vector in $\CC^4$ with components $\xi_{z_k}$,
and the extension with class $\xi$
provides a sheaf $\LLL$ locally free at $z_k$ as an $\OOO_{\G_i}$-module
if and only if $\xi_{z_k}\neq 0$.

Let $s$ be the number of points $z_k$ in which $\LLL$ is
locally free as an $\OOO_{\G_i}$-module. Then $\FFF\simeq  \LLL_{i}(-s\cdot pt)$
and $\FFF'\simeq  \LLL'_{i}$, where $\LLL_{i}=\gamma_i^{-1}\LLL\simeq
\LLL|_{C_i}/\mbox{(torsion)}$, $\LLL'_{i}=\gamma_i^{\prime -1}\LLL\simeq
\LLL|_{C'_i}/\mbox{(torsion)}$, and $\gamma_i$ (resp. $\gamma'_i$) is the
natural inclusion of $C_i$ (resp. $C'_i$) into $\G_i$.
Thus (\ref{triple})
acquires the form
\begin{equation}\label{L-to-C}
0\lra\LLL_{i}(-s\cdot pt)\lra\LLL\lra\LLL'_{i}\lra 0
\end{equation}
Let $d=\deg \LLL_{i}$,
$d'=\deg \LLL'_{i}$. We will call $(d,d')$ the bidegree of $\LLL$.
Then the semistability of $\LLL$
implies $d-s\leq d'$. Reversing the roles of $C_i, C'_i$, we can
represent the same sheaf as an extension
\begin{equation}\label{L-to-C'}
0\lra\LLL'_{i}(-s\cdot pt)\lra\LLL\lra\LLL_{i}\lra 0,
\end{equation}
which implies $d'-s\leq d$. We have $\chi (\LLL)=k-2=d+d'-s+2$ 
and $|d-d'|\leq s\leq 4$. Taking $s=4$, we obtain all the locally
free extensions; the only possible bidegrees
are given by $(d-\frac{k}{2},d'-\frac{k}{2})\in
\{ (-1,3), (0,2), (1,1), (2,0), (3,-1)\}$.
The extremal cases $(-1,3),(3,-1)$ correspond to
$\deg\FFF=\deg\FFF'=\frac{k-4}{2}$,
so all such extensions represent one and the same $S$-equivalence
class of the {\em trivial} extension,
that is the direct sum 
$\OOO_{C_i}(\frac{k-4}{2}pt)\oplus\OOO_{C'_i}(\frac{k-4}{2}pt)$.
For the remaining three bidegrees, the non-isomorphic 
{\em locally-free} extensions provide non-isomorphic stable sheaves.
The locally free extensions are parametrized by the complements
$J^{d,d'}\simeq (\CC^*)^3$
to the coordinate hyperplanes in 
$\PP(\Ext^1(\OOO_{C'_i}(d'pt),\OOO_{C_i}((d-4)pt)))\simeq \PP^3$,
so $J^{d,d'}$ are mapped injectively into $\MMM^k$.
The non-locally-free extensions deform in the corresponding
$\Ext$-groups to the locally free ones, so they lie
in the closures $\bar{J}^{d,d'}$ of the images of $J^{d,d'}$.

(iii) Tensoring by $\OOO_S(q)$ for $q=C_i$ preserves the support and
the property of being torsion-free rank-1 sheaf considered as
a sheaf on its support. Thus it preserves the stability of all
the sheaves from $\MMM^k$ supported on the integral curves.
But it changes the distribution of degrees on the components
of reducible ones. If we denote by $(\tilde{d},\tilde{d'})$
the bidegree of $\LLL(q)$ for $\LLL$ supported on $\Gamma_j$
we have:
$$
(\tilde{d},\tilde{d'})=\left\{\begin{array}{ll}
(d-2,d'+4) & \mbox{if}\ \ j=i \\
(d+2,d') &  \mbox{if}\ \ j\neq i.\end{array}\right.
$$
This immediately implies the formulas for the indeterminacy locus of~$\psi$.
\end{proof}

\begin{remark}
For odd $k=2m+1$, $\MMM^{2m+1}$ is smooth and is birational
to $\MMM^3$. In its turn, $\MMM^3$ is birational to the punctual 
Hilbert scheme $S^{[3]}$ (see \cite{Beau-3}). Then, by
\cite{Hu-0}, $\MMM^{2m+1}$ is deformation equivalent
to $S^{[3]}$.
\end{remark}

\begin{definition}
We will call $\MMM^k$ the degree-$k$ relative compactified  Picard
variety of the linear system $|H|$ and denote it $\PPic^k(|H|)$.
\end{definition}

\section{Local structure of $\PPic^{2m}(|H|)$}
\label{local}

We will use the approach of \cite{O'G-2} to describe the local structure
of the moduli space at a point representing a strictly
semistable sheaf $\FFF$ as a quotient of the versal deformation
of $\FFF$ by $\Aut (\FFF)$.

Let us fix an integer $m$ and consider the relative compactified  Picard
variety $\MMM=\PPic^{2m}(|H|)$.
First we will describe Simpson's construction for $\MMM$.
Let $\LLL\in\MMM$ and
$k\gg 0$ a sufficiently big integer. Then $\LLL (k)$ is generated by
global sections, and denoting by $V$ the vector space $H^0(\LLL (k))$,
we will consider the Grothendieck Quot-scheme $\QQQ $ parametrizing
all the quotients $V\otimes \OOO_X(-k)\onto \LLL'$ 
such that $\chi (\LLL(n))=\chi (\LLL'(n))$ for all $n\in\ZZ$.
Let $\QQQ _c^{ss}\subset \QQQ $ be the open subscheme
parametrizing the semistable pure 1-dimensional sheaves
and $\QQQ _c$ the closure of $\QQQ _c^{ss}$ in $\QQQ $.
There is a natural action of $G=GL(V)$ on $\QQQ $, $\QQQ _c$ and a $G$-linearized
ample invertible sheaf $L$ on $\QQQ $, such that $\QQQ _c^{ss}$ coincides
with the set of $L$-semistable points of the action of $G$ on $\QQQ _c$, and
$\MMM$ is obtained as the Mumford quotient $\QQQ _c//G$.

Let $z\in \QQQ _c^{ss}$ be a point with closed orbit $G\cdot z$,
$[z]$ the corresponding point in $\MMM$,
$\LLL_z$ the quotient sheaf represented by $z$,
and $H$ the stabilizer of $z$; we have $H\simeq\Aut (\LLL_z)$.
Luna's Slice Theorem (\cite{Lu}, \cite{Sim}) affirms that
there exists a $H$-invariant affine subscheme $W\subset \QQQ _c^{ss}$
passing through $z$ such that the map $W//H\lra \QQQ _c//G$ of GIT
quotients is
\'etale. Such a $W$ is called Luna's slice of the action of $G$.
Let $(W,z)$ be the germ of $W$ at $z$ and $\LLLL$ the restriction
of the universal quotient sheaf on $\QQQ \times S$ to $(W,z)\times S$.
By \cite{O'G-2}, Proposition (1.2.3), $((W,z),\LLLL)$ is a versal
deformation of $\LLL_z$.

There is a standard method for constructing a versal deformation
of a sheaf which provides the following proposition:

\begin{proposition}\label{obstruction}
Let $X$ be a smooth projective variety, $\FFF_0$ a coherent
sheaf on $X$. Then there exists a germ of
a nonsingular algebraic variety $(M,0)$ together with a morphism
$\Upsilon:(M,0)\lra (\Ext^2(\FFF_0,\FFF_0),0)$, called
the obstruction map,
such that the following properties are verified:

(i) 
$(\Upsilon^{-1}(0),0)$ is the base of a versal deformation
of $\FFF_0$, that is, 
there exists a coherent sheaf $\FFF$ on
$(\Upsilon^{-1}(0),0)\times X$
such that $((\Upsilon^{-1}(0),0),\FFF)$ is
a versal deformation of $\FFF_0$. The Kodaira-Spencer
map of this deformation provides a natural isomorphism
$KS:T_0M\isoto \Ext^1(\FFF_0,\FFF_0)$.

(ii) Let 
$$
\Upsilon =\displaystyle\sum_{i=1}^\infty\Upsilon_i,\ \ 
\Upsilon_i\in \Hom_{\CC-{\rm lin}} (S_i(T_0M),\Ext^2(\FFF_0,\FFF_0))
$$
be a Taylor expansion of $\Upsilon$. Then $\Upsilon_1=0$
and $\Upsilon_2$ is the composition
\begin{multline*} T_0M\xrightarrow{KS\times KS}
\Ext^1(\FFF_0,\FFF_0)\times\Ext^1(\FFF_0,\FFF_0)
\xrightarrow{(\xi, \xi)\mapsto \xi\smallcup\xi}
\Ext^2(\FFF_0,\FFF_0)
\end{multline*}
where $\xi\smallcup\eta$ denotes the Yoneda product of two
elements of $\Ext^1(\FFF_0,\FFF_0)$.
\end{proposition}

\begin{proof}
The Appendix of Bingener to \cite{BH} provides a scheme of the proof of this statement.
The existence of a formal versal deformation was proven in
\cite{Rim}. By \cite{Art}, the formal versal deformation
is the formal completion of a genuine versal deformation.
The identification of the obstruction $\Upsilon_2$ on the formal level
with the Yoneda pairing was done in \cite{Ar}, \cite{Mu-2}.
See also \cite{HL}, I.2.A.6 and historical comments,
for the case when $\FFF_0$ is simple.
For the construction of $\Upsilon_i$ for all $i$, see Proposition A.1
of \cite{LS}. See also the paper \cite{Lau}, which provides a similar construction in the deformation theory of modules over a $k$-algebra and uncovers its relation to
the Steenrod squares.
\end{proof}

\begin{lemma}
In the situation of Proposition \ref{obstruction},
let us assume in addition that $X$ is a  K3 or abelian surface.
Then the image of $\Upsilon_2$ lies in the codimension-$1$
subspace $\Ext^2(\FFF_0,\FFF_0)_0$, defined as the
kernel of the trace map
$
\Ext^2(\FFF_0,\FFF_0)
\xrightarrow{\Tr}H^2(\OOO_X).
$
\end{lemma}

\begin{proof}
The surjectivity of the trace map on a K3 or abelian surface is proved
in \cite{Ar}, \cite{Mu-2}. The fact that $\Tr\circ \Upsilon_2=0$ 
follows from \cite{HL}, Lemma 10.1.3.
\end{proof}

Let now $z$ be a point of $\QQQ _c^{ss}$ representing
one of the 28 polystable sheaves 
$\LLL_z=\OOO_{C_i}((m-2)pt)\oplus\OOO_{C'_i}((m-2)pt)$.
To shorten the formulas, we will denote 
$\LLL_z=\LLL$, $\LLL_1=\OOO_{C_i}((m-2)pt)$, $\LLL_2=\OOO_{C'_i}((m-2)pt)$,
so that $\LLL=\LLL_1\oplus\LLL_2$. As in the previous section,
denote by $z_1,\dots,z_4$ the intersection points of $C_i$ and $C'_i$.
The orbit of a polystable sheaf is closed in
$\QQQ _c^{ss}$ (see \cite{LP}, 2.9), so
the above local description of $\MMM$ at $[z]$ applies.

We have for $i=1,2$:
$$
\EXT^1_{\OOO_S}(\LLL_i,\LLL_{2-i})\simeq\bigoplus_{q=1}^4\CC_{z_q},\
\EXT^k_{\OOO_S}(\LLL_i,\LLL_{2-i})=0\ \mbox{if}\ k\neq 1,
$$
\begin{multline*}
\EXT^0_{\OOO_S}(\LLL_i,\LLL_{i})=\OOO_C,\ 
\EXT^1_{\OOO_S}(\LLL_i,\LLL_{i})\simeq\OOO_C(-2pt), \\
\EXT^k_{\OOO_S}(\LLL_i,\LLL_{i})=0\ \mbox{if}\ k\not\in \{ 0,1\},
\end{multline*}
where $C=C_i$ for $i=1$ and $C=C'_i$ for $i=2$. Thus
$$
T_zW\simeq \Ext^1(\LLL,\LLL)=\Ext^1(\LLL_1,\LLL_2)
\oplus \Ext^1(\LLL_2,\LLL_1),\ 
\Ext^1(\LLL_i,\LLL_i)=0,$$
$$
\Ext^1(\LLL_i,\LLL_{2-i})=H^0(\EXT^1_{\OOO_S}(\LLL_i,\LLL_{2-i}))\simeq
\bigoplus_{q=1}^4\CC_{z_q},
$$
$$
\Ext^2(\LLL,\LLL)=\bigoplus_{i=1,2}\Ext^2(\LLL_i,\LLL_i),\ \ 
\Ext^2(\LLL_i,\LLL_{2-i})=0,\ \  i=1,2.
$$
By Serre duality (\cite{Mu-2}, Proposition 2.3),
$\Ext^2(\LLL_i,\LLL_i)\xrightarrow{\Tr}H^2(\OOO_S)$ is
an isomorphism, and
$$
\Ext^1(\LLL_i,\LLL_j)\times \Ext^1(\LLL_j,\LLL_i)
\xrightarrow{\rm Yoneda}\Ext^2(\LLL_i,\LLL_i)\xrightarrow{\Tr}H^2(\OOO_S),
$$
where $j=2-i$, is a nondegenerate pairing.
Let us fix once and forever a generator of $H^2(\OOO_S)$,
then denote by $e_i$ its preimage in 
$\Ext^2(\LLL_i,\LLL_i)$. Denote $E_i=\Ext^1(\LLL_i,\LLL_{2-i})$,
$E=\Ext^1(\LLL,\LLL)$.
Our choice of the $e_i$ allows
us to identify $E_{2-i}$ with the dual of $E_i$
in such a way that $E=E_1\oplus E_2$ and $\Upsilon_2$ is given by
\begin{equation}\label{first-ob}
\Upsilon_2:E_1\oplus E_2\lra \CC e_1\oplus\CC e_2,\ \ (\xi_1,\xi_2)
\mapsto \langle\xi_1,\xi_2\rangle (e_1-e_2)
\end{equation}
Thus we have proved:

\begin{lemma}
The first obstruction map $\Upsilon_2$ for the sheaf $\LLL$
is a nondegenerate quadratic form on the $8$-dimensional vector
space $E=\Ext^1(\LLL,\LLL)$ with values in the $1$-dimensional
vector space $\Ext^2(\LLL,\LLL)_0=\CC (e_1-e_2)$,
given by formula (\ref{first-ob}).
\end{lemma}

This implies, in particular, that the base of the versal
deformation $\Upsilon^{-1}(0)$ is at most 7-dimensional.
Further, the stabilizer $H$ of $z$
is just the automorphism group $\Aut (\LLL)=\CC^*\id_{\LLL_1}\times\CC^*\id_{\LLL_2}$
acting on $\Upsilon^{-1}(0)$ via its quotient by the center,
hence with 1-dimensional orbits. As $\dim\MMM =6$, we conclude:

\begin{corollary}\label{quad-sing}
$(\Upsilon^{-1}(0),0)$ is a nondegenerate $7$-dimensional
quadratic singularity with tangent cone $\Upsilon_2^{-1}(0)$.
In particular, $(\Upsilon^{-1}(0),0)$ and $(\Upsilon_2^{-1}(0),0)$
are analytically equivalent.
\end{corollary}

The linearized action of $H$
is given by the following lemma.

\begin{lemma}\label{action}
In the above notation, 
let $g\in H=\Aut (\LLL)$ and \mbox{$g:W\lra W$} the map
given by the group action of $H$ on $W$. Let
$\xi\in T_zW\simeq \Ext^1(\LLL,\LLL)$.
Then $g_*(\xi)=g\smallcup\xi \smallcup g^{-1}$
\end{lemma}

\begin{proof}
See \cite{O'G-2}, (1.4.16), or \cite{Dr}, (7.4.1).
\end{proof}

Thus
the linearized action on $E$ of an element $(\lambda_1,\lambda_2)\in
\Aut (\LLL)$ is given by:
$$
(\lambda_1,\lambda_2)_*:E_1\oplus E_2\lra E_1\oplus E_2, \ \ 
(\xi_1,\xi_2) \mapsto (\lambda_1^{-1}\lambda_2\xi_1,
\lambda_1\lambda_2^{-1}\xi_2).
$$
Passing to the quotient $\mathrm PH=\PAut (\LLL)$ by the center,
we have: \mbox{$\mathrm PH\simeq \CC^*$} via the map
$(\lambda_1,\lambda_2)\mapsto \lambda=\lambda_1^{-1}\lambda_2$,
and for $\lambda\in \mathrm PH$,
$\lambda_*$~acts with weight $\lambda$ on $E_1$ and $\lambda^{-1}$ on $E_2$.
Let us introduce coordinates $x_1,\ldots,x_4$ on $E_1$ in such a
way that the $i$-th coordinate axis is $H^0(\CC_{z_i})\subset E_1$.
Let $y_1,\ldots,y_4$ be the dual coordinates on $E_2$.
We obtain:

\begin{corollary}\label{H-in-coord}
The linearized action of $\mathrm PH$ on $\Upsilon_2^{-1}(0)$
is identified with the action of $\CC^*$ on the nondegenerate
quadraric cone $Q=\{x_1y_1+\ldots +x_4y_4=0\}$ in $E\simeq\CC^8$ given by
\begin{multline}\label{action-of-lambda}
\lambda_*:(x_1,\ldots,x_4,y_1,\ldots,y_4)\ \mapsto\\
(\lambda x_1,\ldots,\lambda x_4,\lambda^{-1}y_1,\ldots,\lambda^{-1}y_4),
\ \ \lambda\in \CC^*.
\end{multline}
\end{corollary}

Now we will use the birational modification $\pi:
\tilde{\MMM}\lra \MMM$ of $\MMM$ constructed by the method of Kirwan \cite{Kir}.
Given a GIT quotient $Z//G$, Kirwan constructs
its partial desingularization in blowing up successively
closed semistable orbits of $G$
until the stability of a point under the action
of $G$ on the blown up variety $\tilde{Z }$ becomes equivalent
to the semistability:\ \ 
$\tilde{Z }^{ss}=\tilde{Z }^{s}$. Equivalently, one may require
that, the projectivized stabilizer of any 
semistable point of $\tilde{Z}$ is finite.
Then {Kirwan's modification} of $Z//G$ is the induced birational
morphism $\pi: \tilde{Z }//G\rar Z//G$.

In our case, we consider just one blowup $\sigma:\tilde{\QQQ }_c\lra \QQQ _c$ 
with center at the union of all the closed
orbits in the strictly semistable locus of $\QQQ _c$. 
The induced map of GIT quotients
$\pi:\tilde{\MMM}=\tilde{\QQQ }_c//G\lra\MMM=\QQQ _c//G$
is a morphism by \cite{Kir}, 3.1, 3.2.

\begin{proposition}\label{kirwan}
In the above notation, the following
assertions hold:

(i)  $\pi :\tilde{\MMM}\lra \MMM$ is Kirwan's 
modification of $\MMM$.

(ii) $\tilde{\MMM}$ is nonsingular and projective, thus $\pi$
is a resolution of singularities of $\MMM$.
The construction of $\tilde{\MMM}$ consists in blowing
up the $28$ singular points
$\zeta_1,\ldots ,\zeta_{28}$ of $\MMM$ taken with their reduced
structure.

(iii) The exceptional divisors
$I_i=\pi ^{-1}(\zeta_i)$ ($i=1,\ldots ,28$) can be identified with
the flag variety $\Fl(0,2;\PP^3)$.

(iv) The normal bundle $\NNN_{I_i/\tilde{\MMM}}$ is isomorphic to
$\OOO_{I_i}(-1)$, the restriction 
of $\OOO_{\PP^{14}}(-1)$ to the flag variety
$\Fl(0,2;\PP^3)$ in its standard embedding into $\PP^{14}$.
\end{proposition}

\begin{proof}
(i) We have to show that one blowup suffices to get a complete
Kirwan's modification.
For a strictly semistable point $z\in \QQQ _c$, $\sigma$ induces
on the \'etale slice $W$ at $z$ the blowup $\sigma_z:\tilde{W}\lra W$
with center $z$, and $\tilde{W}//H$ is an \'etale neighborhood
of the exceptional fiber $\pi^{-1}(\zeta)$ in $\tilde{\QQQ }_c//G$,
where $H$ is the stabilizer of $z$ and $\zeta=[z]$ is the image of $z$
in $\MMM$.
By Corollary \ref{H-in-coord}, $F=\sigma_z^{-1}(z)$ is isomorphic
to the projectivized quadratic cone $\PP Q$ with $\mathrm PH$
acting by formula~(\ref{action-of-lambda}). 
The two projective 3-spaces $\PP E_1$, $\PP E_2$ contained
in $\PP Q$ consist of unstable points, and the stabilizer
of any point of $\PP Q\smallsetminus (\PP E_1\cup\PP E_2)$
in $\mathrm PH$ is $\{\pm 1\}$, so
$\PP Q^{ss}=\PP Q^s=\PP Q\smallsetminus (\PP E_1\cup\PP E_2)$.
As all the semistable points of $\PP Q$ are stable, $\pi$
is Kirwan's modification at~$\zeta$.
Remark also that the strictly semistable
points of $\QQQ _c$ (or $W$) with {\em non-closed} orbits become unstable
when lifted to $\tilde{\QQQ }_c$ (resp.$\tilde{W}$).

(ii) 
The blowup $\tilde{W}$ at $z$ is nonsingular over $z$ since,
by Corollary \ref{quad-sing}, $(W,z)$ is a nondegenerate
quadratic singularity. As the stabilizer in $\mathrm PH$
of all the semistable
points of $\sigma^{-1}(z)$ is constant, equal to $\{\pm 1\}$,
the quotient $\tilde{W}//H$ is nonsingular at every
point of $F=\pi^{-1}(z)$ by Luna's slice theorem. By Lemma 3.11 of \cite{Kir},
$\pi$ is the blowup of the reduced point $\zeta=[z]$.

(iii) The exceptional fiber $I=\pi^{-1}(\zeta)$ is isomorphic
to the quotient $\PP Q//\CC^*$ by the action
(\ref{action-of-lambda}). The algebra of invariants of
this action is generated by the quadratic monomials
$u_{ij}=x_iy_j$, and the generating relations are of two types:
one linear, $u_{11}+\ldots +u_{44}=0$, and the quadratic ones
$u_{ij}u_{kl}=u_{kj}u_{il}$. The quadratic relations define
the standard Segre embedding of $\PP^3\times\PP^3$ in $\PP^{15}$,
and the linear one cuts out the incidence variety: if we identify
the second factor $\PP^3$ with $\PP^{3\dual}$,
parametrizing the hyperplanes $h$ in the
first factor $\PP^3$, then $I=\{ (p,h)\in\PP^3\times\PP^{3\dual}
\mid p\in h\}$. This is just the flag
variety $\Fl(0,2;\PP^3)$.

(iv) Let $A$ denote the algebra of regular functions on $W$,
so that $W=\Spec A$.
Let $\mathfrak M=\mathfrak M_z\subset A$ 
be the maximal ideal of $z\in W$.
As any representation of $\mathrm PH$ is completely reducible, $\mathfrak M$
contains a subrepresentation $V$ of $\mathrm PH$ that projects down isomorphically
and equivariantly onto $\mathfrak M/\mathfrak M^2=T_z^*W$ via the differential
$d:\mathfrak M\rar\mathfrak M/\mathfrak M^2$. The map $d|_V$ extends
to a $\mathrm PH$-equivariant epimorphism of $\CC$-algebras $S_\ccd V\rar A$
giving rise to a $\mathrm PH$-equivariant morphism $W\rar T_zW$.
Its image is a hypersurface $W_1$, and it is \'etale at $z$ if considered as
a morphism $W\rar W_1$.
In shrinking $W$, we can assume that
this morphism is everywhere \'etale.

We can choose
functions $x_i,y_i\in V\simeq T^*_zW$ on which $\mathrm PH$ acts according to
formula (\ref{action-of-lambda}). Then $W_1$ is defined by the equation $F=0$,
where $F\in\CC[x_i,y_i]$ is the sum of homogeneous components
of even degree, $F=F_2+\ldots +F_{2r}$, which are $\mathrm PH$-invariant
and such that $F_2=x_1y_1+\ldots +x_4y_4$. We can write $W_1=\Spec A_1$,
where $A_1=\CC[x_i,y_i]/(F)$. Let $U=W//\mathrm PH$, $U_1=W_1//\mathrm PH$.
We have $U=\Spec B$, $U_1=\Spec B_1$, where $B=A^{\mathrm PH}$,
$B_1=A_1^{\mathrm PH}$, and the \'etale morphism $W\rar W_1$
descends to the quotients as an \'etale morphism $U\rar U_1$.
Let $\zeta=[z]$ denote the image of $z$ in $U$ or in $U_1$,
and $\mathfrak m=\mathfrak m_\zeta$ the maximal ideal of $\zeta$
in either one of the rings $B, B_1$. 
The blowup $\tilde U$ of $U$ at $\zeta$ can be given by
$\tilde U=\Proj_B(\bigoplus\limits_{k\geq 0}\mathfrak m^k)$, and the exceptional
divisor $I=\Proj_\CC(\bigoplus\limits_{k\geq 0}\mathfrak m^k/\mathfrak m^{k+1})$,
its normal bundle being $\OOO_I(-1)$, the dual of the Grothendieck
tautological sheaf $\OOO_I(1)$ on the latter Proj.
As $U$ and $U_1$ are locally isomorphic at $z$ in the \'etale topology,
$I$ and its normal bundle do not depend on whether $\mathfrak m$ is considered
in $B$ or in $B_1$. So the wanted normal bundle
$\NNN_{I_i/\tilde{\MMM}}$ can be computed as the
normal bundle to the blowup of $\zeta $ in $U_1$.

Choosing $u_{ij}=x_iy_j$ as the generators of $B_1=A_1^{\mathrm PH}$,
we represent $U_1$ as a hypersurface in the cone $\mathfrak C=\{(u_{ij})\in\CC^{16}
\mid u_{ij}u_{kl}=u_{kj}u_{il},\ \  1\leq i,j,k,l\leq 4\}$,
defined by the equation $f=0$, where $f$ has a decomposition into homogeneous
components of the form $f=f_1+f_2+\ldots+f_r$, $f_1=u_{11}+\ldots +u_{44}$. 
The proper transform $\tilde U_1$ of the hypersurface $f=0$ in the blowup $\tilde{\mathfrak C}$
of $\mathfrak C$ at $\zeta$ meets the exceptional divisor
$E\simeq \PP^3\times\PP^3$ transversely along the
flag variety $I\subset E$. Hence $\NNN_{I/\tilde U_1}\simeq
\left(\NNN_{E/\tilde{\mathfrak C}}\right)|_I$. But the latter normal bundle is
just the restriction of $\OOO_{\PP^{15}}(-1)$, and we are done.
\end{proof}

\begin{remark}
Our argument in part (iv) is a kind of ``equivariant deformation to the normal
cone'', compare to Sect. 5 of \cite{LS}.
\end{remark}

The exceptional divisor $I_j$ over any of the points $\zeta_j$ has two
distinct projections to $\PP^3$ which are $\PP^2$-bundles, and
which we will refer to as {\em rulings} of $I_j$.
By Proposition \ref{kirwan} (iv),
the normal bundle to $I_j$ restricts as $\OOO (-1)$
to the fibers $\PP^2$ of each ruling. 
By Moishezon's contractibility criterion \cite{Mo},
both projections of $I_j$ to $\PP^3$ can be extended to
a morphism $f:\tilde{\MMM}\lra Y$ such that $Y$ is a smooth
compact complex 6-dimensional manifold, not necessarily projective.
Applying this argument successively to different $j=1,\ldots , 28$,
we obtain:

\begin{corollary}\label{2-to-28}
There are $2^{28}$ distinct bimeromorphic morphisms
\mbox{$f:\tilde{\MMM}\lra Y$} onto smooth, compact, complex,
not necessarily projective 6-dimensional manifolds $Y$
which contract each one of the divisors $I_j$ onto a
projective $3$-space $f(I_j)\simeq\PP^3$.
For any of these morphisms $f$, 
Kirwan's desingularization $\pi:\tilde{\MMM}\lra \MMM$ factors through $f$,
$\pi =g\circ f$, so that $g$ is a small contraction, that is,
a contraction without exceptional divisors. Moreover,
the symplectic form $\alp$ on the nonsingular locus $\MMM^*$
of $\MMM$ lifts to a global symplectic form $\alp_Y$ on $Y$,
and hence $Y$ is a holomorphically symplectic manifold.
\end{corollary}

\begin{proof}
The small contraction map $g$ induces an isomorphism
$g:g^{-1}(\MMM^*)\isoto \MMM^*$, so $g^*(\alp)$ is a symplectic
form on $g^{-1}(\MMM^*)$. It extends to a regular 2-form
$\alp_Y$ on all of $Y$ by Riemann--Hartogs extension
theorem since the complement $Y\smallsetminus g^{-1}(\MMM^*)$ is a union
of $\PP^3$'s and thus is of codimension $>1$. Finally, $\alp_Y$
is nondegenerate, and hence is a symplectic form.
Indeed, the degeneracy locus of $\alp_Y$ is
nothing else but the zero locus of $\alp_Y^{\wedge 3}\in H^0(Y,\Omega^6_Y)$. 
The zero locus of a section of an invertible sheaf, if nonempty, is
either $Y$ itself, or a divisor in $Y$, but we know that
$\alp_Y$ is nondegenerate on an open set whose complement
contains no divisors, so $\alp_Y$ is everywhere nondegenerate.
\end{proof}

In fact, there are projective varieties among the complex manifolds
$Y$ from Corollary \ref{2-to-28}. One of them is given by the next
proposition.

\begin{proposition}
Let $H_\eps=H+\eps\sum_{i=1}^{28}(C_i-C'_i)$.
Then there exists a sufficiently small $\eps_0>0$ such that for any 
$\eps\in\QQ\cap ]0,\eps_0[$, the following assertions hold:

(i) $H_\eps$ is an ample $\QQ$-divisor on $S$.

(ii) The (semi-) stability of a sheaf with Mukai vector $v=(0,H,2m-2)$ 
with respect to $H_\eps$ does not depend on $\eps$,
and every $H_\eps$-semistable sheaf with Mukai vector $v$ is stable.

(iii) The moduli space $Y=M_S^{H_\eps,ss}(v)=M_S^{H_\eps,s}(v)$
is an irreducible symplectic manifold which does not depend on $\eps$.

(iv) The natural map $g:Y\rar \MMM$ is a small resolution of singularities
such that $g^{-1}(\zeta_j)\simeq\PP^3$ ($j=1,\ldots,28$).
\end{proposition}

\begin{proof}
(i) follows from the openness of the cone of ample classes in $\Pic S\otimes\RR$.
For (ii), remind that the (semi)-stability of a sheaf supported on an integral
curve does not depend on polarization. So we have only to examine the sheaves
supported on the reducible curves $\G_i$. This is similar to
the proof of Proposition \ref{fibers-of-f}, (iii): any $H_\eps$-semistable sheaf
which is rank-1 and torsion-free as a sheaf on its support $\G_i=C+C'$
is given by extensions (\ref{L-to-C}) and (\ref{L-to-C'}) such that
\begin{multline*}
(1-3\eps)(d-s)\leq (1+3\eps)d',\ \ 
(1+3\eps)(d'-s)\leq (1-3\eps)d,\\ s=0,\ldots,4,\ \ d+d'=2m. 
\end{multline*}
For all sufficiently small
$\eps>0$, the solutions of these inequalities are
the same triples of integers $d,d',s$ as in the proof of Proposition \ref{fibers-of-f}, (iii), 
except for $d=d'=m$, $s=0$ which does not satisfy the second inequality. 
For all the solutions, the inequalities are strict, hence
there are no strictly semistable sheaves. This ends the proof of (ii).
The assertion (iii) follows by \cite{HL}, 6.2.5. To prove (iv), remark, that
by the above argument, any $H_\eps$-semistable sheaf is also $H$-semistable, so
there is a natural morphism $g:Y\rar \MMM$. Further,
all the nontrivial extensions
$$
0\lra \OOO_C((m-2)pt)\lra \FFF\lra \OOO_{C'}((m-2)pt)\lra 0
$$
provide $H_\eps$-stable sheaves with the same image $[\OOO_C((m-2)pt)\oplus
\OOO_{C'}((m-2)pt)]=\zeta_i$ in $\MMM$,
and two such sheaves are isomorphic if and only if they correspond
to proportional extension classes.
Thus $g^{-1}(\zeta_i)=\PP\Ext^1(\OOO_{C'}((m-2)pt),
\OOO_C((m-2)pt))\simeq\PP^3$.
\end{proof}

\section{The relative compactified  Prymian $\PPrym^k(\CCC,\tau)$}
\label{prymian}

The Galois involution $\tau$ of the double covering
$\rho:S\rar X$ is $H$-linear and induces an involution on
$|H|\simeq\PP^3$, whose fixed locus consists of two components:
a point and a plane. The plane
parametrizes the curves of the form $\rho^{-1}\mu^{-1}(\ell)$,
where $\ell$ runs over the lines in $\PP^2$, thus this plane
is identified with the dual of the $\PP^2$ which is the image of
$\mu$. We will denote it $\PP^{2\dual}$. The other component
of the fixed locus, a point, corresponds
to the ramification curve $\Delta\in |H|$.
Let $\phi:\CCC\rar\PP^{2\dual}$ be the linear subsystem of
$\tau$-invariant curves parametrized by $\PP^{2\dual}$. A generic
$t\in \PP^{2\dual}$ represents a line $\ell=\ell_t$ which is not
tangent to $B_0$, neither to $\bar{\Delta}_0:=\mu(\Delta_0)$.
The corresponding curve $C_t=\phi^{-1}(t)=\rho^{-1}\mu^{-1}(\ell_t)$
is a smooth genus-3 curve, and $E_t=C_t/\tau$ is elliptic;
the double cover $\rho_t=\rho|_{C_t}:C_t\rar E_t$ is branched at 4 points
of the intersection $\Delta_0\cap E_t$.

\begin{definition}
Let $\eta:C\rar D$ be a double covering map of integral projective curves
and $\tau$ the Galois involution of $\eta$.
Then $\tau$ acts as a linear involution on the generalized Jacobian
$J(C)$, and the Prym variety $\Prym (C,\tau)$ is defined as $\im (\id-\tau)=
[\ker (\id+\tau)]^\circ$, where $G^\circ$ denotes the connected
component of the neutral element in a subgroup $G$ of $J(C)$.
\end{definition}

If $C$ is smooth, then $J(C)$ and $\Prym (C,\tau)$ are abelian varieties,
but for singular curves, they can be extensions of abelian varieties by a number
of copies of $\CC^*$ or $\CC$.

\begin{lemma}\label{connected-kernel}
Let $C$ be a smooth genus-$3$ curve with an involution $\tau$
such that $D=C/\tau$ is an elliptic curve. 
Then $\ker (\id+\tau)$ has only one connected component in $J(C)$, and
the restriction of the principal polarization from $J(C)$ to
$\Prym (C,\tau)=\ker (\id+\tau)$ is a polarization of
type $(1,2)$.
\end{lemma}

\begin{proof}
It is well-known that under the hypotheses of the lemma,
$P=\Prym (C,\tau)$ has a polarization of type $(1,2)$,
see \cite{B}.
A very explicit proof of the fact that this polarization is the
restriction of the standard principal polarization of the Jacobian
$J=J(C)$ is given in the paper \cite{P}, in which the author identifies
the intersection $\Theta_a\cap P$, where $\Theta_a=a+\Theta$ is
an appropriate translate of the theta-divisor $\Theta\subset J$, as
a genus-3 curve $C^\dual$ obtained by a bigonal construction from $C$.

Let $\eta:C\rar D$ be the natural double covering map.
As $C$ is not hyperelliptic, $\eta^*:J(D)\rar J$ is injective.
Let $E=\eta^*(JD)\subset J$. Then $E+P=J$, and $K=E\cap P\subset J_{(2)}$,
where $J_{(n)}$ denotes the $n$-torsion subgroup of $J$. It is obvious that
$\ker (\id+\tau)=\bigcup_{z\in E_{(2)}}(z+P)$, so $\ker (\id+\tau)$ is
connected if and only if $K=E_{(2)}$. By \cite{BM}, 7.6 and 7.10,
$K=\ker \lambda_{\Xi^1}=\ker \lambda_{\Xi^2}$, where $\Xi^1=\Theta|_E$,
$\Xi^2=\Theta|_P$, and $\lambda_\Xi$ denotes the polarization map
associated to an ample divisor $\Xi$ on an abelian variety $A$. It is
defined by the formula
$$
\lambda_\Xi:A\lra \hat A=\Pic^0(A),\ \ a\mapsto \Cl[\Xi_a-\Xi].
$$
Since we already know that $\Xi^2$ is a polarization of type $(1,2)$,
we have $\#(\ker \lambda_{\Xi^2})=4$, hence $K=E_{(2)}$ and we are done.
\end{proof}

The lemma allows us to define $\Prym(C,\tau)$ as the fixed locus of the
involution $\kappa=\tau\circ\iota$, where $\iota:J(C)\rar J(C)$ is defined
by $[\LLL]\mapsto[\LLL^{-1}]$. Now we will relativize the construction of
$\kappa$ in the linear system $|H|$.

Let $\MMM^k=\PPic^k(|H|)$ be as in the previous sections. 
First, let $k=2m$ be even.
The naive extension of $\iota$ to the sheaves that are not invertible on $C=\Supp\LLL$
is $[\LLL]\mapsto[\HOM_{\OOO_C}(\LLL,\OOO_C(mH))]$. But this does
not commute with base change. The proper definition is
$$
\iota:\MMM^{2m}\lra\MMM^{2m},\ \ [\LLL]\mapsto[\EXT^1_{\OOO_S}(\LLL,\OOO_S((m-1)H))].
$$
This duality functor for pure 1-dimensional sheaves was used by Maruyama in \cite{Maru}, Proposition 2.9, over $\PP^2$, but it can be applied on any smooth surface.
The fixed locus $\Fix(\kappa)$ of $\kappa=\tau\circ\iota$ has one connected component of dimension 4, parametrizing sheaves with supports $C_t$, $t\in\PP^{2\dual}$, and $2^6=64$
isolated points representing the invertible sheaves $\LLL$ on $\Delta$
such that $\LLL^2\simeq\OOO_S(mH)|_\Delta$.

To define $\iota$ for $k=2m+1$, we need to fix a class $c\in\Pic(S)$ of degree 2,
that is such that $(c\cdot H)=2$. Then we define
$$
\iota=\iota_c:\MMM^{2m+1}\dasharrow\MMM^{2m+1},\ \ [\LLL]\mapsto[\EXT^1_{\OOO_S}(\LLL,\OOO_S((m-1)H-c))].
$$
This is a rational involution whose indeterminacy locus consists
of \mbox{$H$-stable} sheaves $\LLL\in \MMM^{2m+1}$ such that $\LLL\otimes\OOO_S(-c)$ is
unstable. One can choose for $c$ one of the 56 conics $C_i,C'_i$. For example, if $c=C'_i$,
then the indeteminacy locus of $\iota$ coincides with $\Indet\psi$ as
given by formula (\ref{indetpsi}). Thus
$\kappa=\tau\circ\iota$ is a rational involution in this case and we define
$\Fix(\kappa)$ as the closure of the fixed point set of the restriction of $\kappa$
to its regular locus. It also has one 4-dimensional and 64 zero-dimensional components.

\begin{definition}\label{def-prym}
The relative compactified  Prymian $\PPrym^{k,\kappa}(\CCC,\tau)$, or simply
$\PPrym^{k}(\CCC,\tau)$, is the $4$-dimensional component of $\Fix(\kappa)$
in $\MMM^k$.
\end{definition}

We will study in more detail the variety $\PPrym^{k}(\CCC,\tau)$
for even $k=2m$, which we will denote by ${\PPPP}^{2m}$, or simply
${\PPPP}$ when there is no risk of ambiguity. Remark that
${\PPPP}^{2m}\simeq{\PPPP}^{2m+4}$ via the map $\FFF\mapsto
\FFF (H)$, so that there are at most two different varieties
${\PPPP}^{2m}$:\ \ \ ${\PPPP}^{0}$ and ${\PPPP}^{2}$. We~ignore if
they are really non-isomorphic, or even non-birational.

\begin{theorem}\label{struct-prym}
Let ${\PPPP}=\PPrym^{2m}(\CCC,\tau)$ with $m\in\ZZ$. Identifying, as above,
the $2$-dimensional linear subsystem of $\tau$-invariant
curves in $|H|$ with $\PP^{2\dual}$, let $f_{\PPPP}=f^{2m}_{\PPPP}$ denote
the map ${\PPPP}\rar \PP^{2\dual}$ sending each sheaf to its support. 
Let $C_t=\phi^{-1}(t)$, $E_t=C_t/\tau$, and $\rho_t=\rho|_{C_t}:C_t\rar E_t$, where
$\phi:\CCC\rar \PP^{2\dual}$ is the natural map and $t\in\PP^{2\dual}$.

Then the following assertions hold:

(i) ${\PPPP}$ is nonsingular out of the $28$ points $\zeta_i=[\LLL_i]$
representing the $S$-equivalence classes of the sheaves $\LLL_i=\OOO_{C_i}((m-2)pt)\oplus
\OOO_{C'_i}((m-2)pt)$, $i=1,\ldots,28$. The singularities $({\PPPP},\zeta_i)$ 
are analytically equivalent to $(\CC^4/\{\pm 1\},0)$.

(ii) ${\PPPP}$ is a symplectic $V$-manifold, and $f_{\PPPP}$ is a Lagrangian fibration on
it. The generic fiber $f_{\PPPP}^{-1}(t)$ is the $(1,2)$-polarized
Prym surface $\Prym(C_t,\tau)$ of the double covering $\rho_t:C_t\rar E_t$.
\end{theorem}

\begin{proof}
(i) It is obvious that the fixed point set of any biregular involution
on a smooth variety is also smooth. The sheaves $\LLL_i$ are invariant
under $\tau$ and $\iota$, hence also under $\kappa$.
So $\zeta_i\in {\PPPP}$, and we only have to determine the analytic type
of the singularity at $\zeta_i$. To this end, we will write out the action
of $\kappa$ on the tangent cone of $\MMM^{2m}$ at $\zeta_i$.

Let us change, for convenience,
the notation, so that $C_+=C_i$, $C_-=C'_i$, $C=\G_i=C_+\cup C_-$, $\LLL_\pm=
\OOO_{C_\pm}((m-2)pt)$, $\LLL=\LLL_+\oplus\LLL_-$, $\zeta=\zeta_i=[\LLL]$.
As $\tau$ leaves invariant
both curves $C_\pm$, it has two fixed points on each of them, which we will
denote by $\lambda_{1\pm},\lambda_{2\pm}$. We can choose homogeneous
coordinates $(x_{0\pm},x_{1\pm})$ on $C_\pm\simeq\PP^1$ in such a way that
$\lambda_{1\pm}=(0:1)$, $\lambda_{2\pm}=(1:0)$, and $\tau$ is given by
$\tau :(x_{0\pm},x_{1\pm})
\mapsto (x_{0\pm},-x_{1\pm})$.
As the
cross-ratio of 4 points of intersection of two conics is the same on
both of them, we can adjust the choice of the above coordinates
in such a way that the $4$ points $z_1,\ldots,z_4$ of $C_+\cap C_-$ have equal coordinates on both curves, and we will number them in such a way that 
they are permuted by $\tau$ in pairs
$z_1\leftrightarrow z_2$, $z_3\leftrightarrow z_4$.

The 4-dimensional vector space $F=\Ext^1(\LLL_+,\LLL_-)$ parametrizes the extensions
$0\rar\LLL_-\rar\FFF\rar\LLL_+\rar 0$. Let $x_i$, resp. $y_i$ be the coordinates
on $F$, resp. $F^\dual=\Ext^1(\LLL_-,\LLL_+)$ obtained in the same way
as those used in Corollary \ref{H-in-coord}.
The choice of $x_i$ made in Section \ref{local} is not unique, it depends on the
choice of a basis in each of the 1-dimensional stalks $\CC_{z_i}$ of the sheaf
$\EXT^1(\LLL_+,\LLL_-)$. Now we will make this choice
more precise. Let $s$ be the number of the points $z_i$ in which $\FFF$ is locally
free. Then $\FFF$ is the result of gluing of the sheaves $\LLL_-(s\cdot pt)$ and $\LLL_+$.
The gluing
consists in the identification of the fibers at $z_i$
via isomorphisms $\phi_i:\LLL_{+,z_i}=\LLL_+\otimes\CC_{z_i}
\isoto \LLL_-(s\cdot pt)_{,z_i}=\LLL_-(s\cdot pt)\otimes\CC_{z_i}$
for those $i$, for which the $\CC_{z_i}$-component of the extension
class is non-zero. Let us denote the resulting
sheaf $\FFF$ by $\LLL_-(s\cdot pt)\#_{(\phi_i)}\LLL_+$.

Consider the case $s=4$. Let us fix some isomorphisms
$\LLL_-(4\cdot pt)\simeq \OOO_{\PP^1}(m+2)$ and $\LLL_+\simeq \OOO_{\PP^1}(m-2)$.
Now, fix $e_+=x_{0+}^{m-2}$, resp. $e_-=x_{0-}^{m+2}$ as a trivializing section
of $\LLL_+$, resp. $\LLL_-(4\cdot pt)$ over an open set containing all the points $z_i$.
Define the four isomorphisms $\phi_i$ as above by $e_{+,z_i}\mapsto e_{-,z_i}$.
Finally, we fix the choice of $(x_i)$ by the condition that $(x_i)$ are the coordinates
of the extension class of the sheaf $$\FFF_{(x_1,x_2,x_3,x_4)}=
\LLL_-(s\cdot pt)\#_{(x_i\phi_i)}\LLL_+$$
whenever $x_i\neq 0$ for all $i=1,\ldots, 4$. This determines also the coordinates
$y_i$, dual to $x_i$.

The action of $\tau$ lifts to $\LLL_\pm$ in such a way that it
preserves $e_-$, $e_+$. Further, $\tau$ interchanges $z_1$ with $z_2$,
$z_3$ with $z_4$, hence $\tau^*(\FFF_{(x_1,x_2,x_3,x_4)})\simeq
\FFF_{(x_2,x_1,x_4,x_3)}$. From here we deduce the action of $\tau$
on $E=\Ext^1(\LLL,\LLL)=F\oplus F^\dual$:
$$
\tau : (x_1,x_2,x_3,x_4,y_1,y_2,y_3,y_4)\mapsto
(x_2,x_1,x_4,x_3,y_2,y_1,y_4,y_3).
$$
As  $\iota$ interchanges $x_i$ with $y_i$, we obtain
$$
\kappa :(x_1,x_2,x_3,x_4,y_1,y_2,y_3,y_4)\mapsto
(y_2,y_1,y_4,y_3,x_2,x_1,x_4,x_3).
$$
The tangent cone to $\MMM^{2m}$ is obtained by taking the
quotient of the quadric $\sum x_iy_i=0$ by $\CC^*$ (see the proof
of Proposition \ref{kirwan}, (iii)). The quotient is identified with
the cone over the hyperplane section $\sum u_{ii}=0$ of
the Segre variety, given by the parametrization $u_{ij}=x_iy_j$ in $\PP^{15}$. 
As we have already noticed, this hyperplane section 
is the flag variety $\Fl(0,2;\PP^3)$ embedded in $\PP^{14}$. Restricting further to the
fixed locus of $\kappa$ is equivalent to intersecting with 6 hyperplanes
$$
u_{11}=u_{22},\ u_{33}=u_{44},\ u_{13}=u_{42},\ u_{14}=u_{32},\ u_{23}=u_{41},\ 
u_{24}=u_{31}.
$$
These equations cut out the Veronese image of $\PP^3$ in $\PP^9$.
Thus the tangent cone to $\MMM^{2m}$ is the cone over the the Veronese image of $\PP^3$,
or in other words, the quotient $\CC^4/\pm 1$. This ends the proof of (i).

The assertions of (ii) follow from (i), Lemma
\ref{connected-kernel} and \cite{Mar-1}, \hyphenpenalty=1000\ Section~6.
\end{proof}

We can use some of the settings of the above proof to determine
the fiber of $f_{\PPPP}$ over a point $t\in\PP^{2\dual}$ representing
a reducible quartic $C_t$.

\begin{lemma}\label{prym-cc}
Let, in the notation of Theorem \ref{struct-prym},
$t\in\PP^{2\dual}$ be a point representing a reducible quartic
$C=C_+\cup C_-$, and $P^{2m}=(f_{\PPPP}^{2m})^{-1}(t)$. Then $P=P^{2m}$ is an irreducible
projective surface having a stratification $P=P_0\sqcup P_1\sqcup P_2$
such that $P_0\simeq\CC^*\times\CC^*$,
$P_1\simeq \CC^*\sqcup\CC^*$ and $P_2$ is a single point. The isomorphism
class of $P$ does not depend on $m$.
\end{lemma}

\begin{proof}
We choose the coordinates $(x_{0\pm},x_{1\pm})$ on $C_\pm$ as in the
proof of Theorem \ref{struct-prym}. Let$z_\pm=x_{1\pm}/x_{0\pm}$ denote the
associated affine coordinate on $C_\pm\setminus\{\lambda_{1\pm}\}$. Then
the $4$ points of $C_+\cap C_-$ have the same values
of the coordinates $z_\pm$: \ $z_{i+}=z_{i-}=z_i$ for $i=1,\ldots,4$.
Moreover, $z_2=-z_1$, $z_4=-z_3$, for the involution $\tau$ acts by
$z_+\mapsto -z_+$, $z_-\mapsto -z_-$. 

Consider first the case $m=0$.
Any invertible sheaf of degree 0 on $C$ can be represented
as the result of gluing $\OOO_{C_+}(a\, pt)$ with $\OOO_{C_-}(-a\, pt)$
at the $4$ points of $C_+\cap C_-$. Let us fix the convention that
the sheaf $\OOO_{C_\pm}(n\, pt)$ is trivialized by the rational
section $e_\pm=x_{0\pm}^{n}$. Denote
the result of the above gluing via the maps $e_{+,z_i}\mapsto \lambda_ie_{-,z_i}$
as $\FFF(a;\lambda_1,\ldots,\lambda_4)$ or
$\OOO_{C_-}(-a\, pt)\connsum{\lambda_1,\ldots,\lambda_4}\OOO_{C_+}(a\, pt)$.
We have  $\FFF(a;(\lambda_i))\simeq \FFF(a';(\lambda_i'))$ if and only if
$a'=a$ and $\lambda_i'=c\lambda_i$ ($i=1,\dots,4$) for some $c\in\CC^*$.
Since $\tau^*$ preserves the value of $a$ and $\iota:\FFF\mapsto \FFF^\dual$
changes $a$ to $-a$, we have: $\FFF(a;(\lambda_i))\in P$
 \ $\Rightarrow$ \  $a=0$. Let $P_0$
be the open subset in $P$ that parametrizes the locally free sheaves.
To determine $P_0$, we compute the action of $\kappa=\iota\circ\tau^*$ on
$\FFF=\OOO_{C_-}\connsum{\lambda_1,\ldots,\lambda_4}\OOO_{C_+}$:
$$
\FFF\sendsto{\tau^*}
\OOO_{C_-}\connsum{\lambda_2,\lambda_1,\lambda_4,\lambda_3}\OOO_{C_+}
\sendsto{\iota}\OOO_{C_-}\connsum{
\frac{1}{\lambda_2},\frac{1}{\lambda_1},\frac{1}{\lambda_4},\frac{1}{\lambda_3}}
\OOO_{C_+},
$$
thus $\FFF\in P_0$ \ $\Leftrightarrow$ \ $\rk\left(
\begin{array}{cccc}
\lambda_1&\lambda_2&\lambda_3&\lambda_4\\
\lambda_2^{-1}&\lambda_1^{-1}&\lambda_4^{-1}&\lambda_3^{-1}
\end{array}\right)=1$, and $P_0$ is the quotient of the subtorus
of $(\CC^*)^4$ with equation $\lambda_1\lambda_2=\lambda_3\lambda_4$
by the diagonal action of $\CC^*$. Hence $P_0\simeq
\CC^*\times\CC^*$.

\hyphenpenalty=10000
We define the next stratum, $P_1$, as the locus of the sheaves in $P$
that are invertible in at least one of the points $z_i$, but are
not invertible in all of them. If $\FFF\in P_1$, then either
$\FFF\simeq \FFF'(a;\lambda_3,\lambda_4)=\OOO_{C_-}((-a-1)pt)\connsum{\, \cdot\, ,\, \cdot\, ,\lambda_3,\lambda_4}\mbox{$\OOO_{C_+}((a-1)pt)$}$, or
$\FFF\simeq \FFF''(a;\lambda_1,\lambda_2)=\mbox{$\OOO_{C_-}((-a-1)pt)$}\connsum{\lambda_1,\lambda_2,\, \cdot\, ,\, \cdot\,}\OOO_{C_+}((a-1)pt)$,
where we put a dot on the $i$-th place
to indicate that the gluing in $z_i$ is not effectuated, that is,
$\FFF_{z_i}=\OOO_{C_-}((-a-1)pt)_{z_i}\oplus\OOO_{C_+}((a-1)pt)_{z_i}$.

To determine the dual of such a sheaf, represent it as the direct
image $\sigma_*(\LLL)$, where $\sigma:\tilde S\rar S$ is the
blowup of $S$ at the two points in which $\FFF$ is not locally free,
and $\LLL$ is supported on the proper transform $C'$ of $C$
and is invertible as a $\OOO_{C'}$-module. Then, by the relative duality
for $\sigma$ (see \cite{Ha}, p. 210),
$$
\FFF^\dual\simeq \sigma_*(
\EXT^1_{\OOO_{\tilde S}}(\LLL,\OOO_{\tilde S}(-\sigma^*(C))\otimes\omega_{\tilde S/S}))
\simeq \sigma_*(\LLL^\dual(-E\cdot C')),
$$
where $E$ is the union of two $(-1)$-curves which form the
exceptional locus of $\sigma$. Let, for example, $\FFF
= \FFF'(a;\lambda_3,\lambda_4)$. Then $\LLL$ is
the result of gluing $\OOO_{C'_-}((-a-1)pt)\connsum{ \lambda_3,\lambda_4}\mbox{$\OOO_{C'_+}((a-1)pt)$}$
at the two points of $C_+'\cap C_-'$, where $C_\pm'$ is the proper
transform of $C_\pm$, and $\LLL^\dual=\LLL^{-1}\simeq
\OOO_{C'_-}((a+1)pt)\connsum{ \lambda_3^{-1},\lambda_4^{-1}}\mbox{$\OOO_{C'_+}((1-a)pt)$}$.
Thus we obtain:
\begin{multline*}
\FFF'(a;\lambda_3,\lambda_4)^\dual=\sigma_*\Bigl(
\OOO_{C'_-}((a-1)  pt)\connsum{ \lambda_3^{-1},
\lambda_4^{-1}}\mbox{$\OOO_{C'_+}((-a-1) pt)$}\Bigr)\simeq\\ 
\FFF'(-a;\lambda_3^{-1},\lambda_4^{-1}).
\end{multline*}
It is much easier to determine the action of $\tau$:
obviously, $\tau^*\bigl(\FFF'(a;\lambda_3,\lambda_4)\bigr)\simeq
\FFF'(a;\lambda_4,\lambda_3)$. We conclude that 
$\FFF'(a;\lambda_3,\lambda_4)\in P_1$
if and only if $a=0$, and the sheaves $\FFF'(0;\lambda_3,\lambda_4)$
fill a component of $P_1$, isomorphic to $\CC^*$. The other copy
of $\CC^*$ contained in $P_1$ is formed by the sheaves
$\FFF''(0;\lambda_1,\lambda_2)$.

Finally, define $P_2$ as the locus of sheaves which are noninvertible
in all the 4 points $z_i$. By Theorem \ref{struct-prym}, 
${\PPPP}$ contains only one such sheaf with support $C$: \ $\FFF=
\OOO_{C_-}(-2pt)\oplus\mbox{$\OOO_{C_+}(-2pt)$}$. This ends the
proof for $m=0$. As $\PPPP^{2m}\simeq\PPPP^{2m+4}$, it remains to consider
the case $m=1$. An isomorphism $P^0\isoto P^2$ can be given by
$\FFF\mapsto\FFF\otimes\theta$, where $\theta$ is a $\tau$-invariant
theta-characteristic on $C$. One easily verifies that there are two such
theta-characteristics: $\theta=
\OOO_{C_-}(pt)\connsum{1,1,\eps,\eps}\OOO_{C_+}(pt)$, where $\eps=\pm 1$.

\end{proof}

\section{Compactified Prymians of integral curves}

We will use the notation of the previous section.
Thus ${\PPPP}$ will denote $\PPrym^{2m}(\CCC,\tau)$, and
$f_{\PPPP}$ or $f^{2m}_{\PPPP}$
the map ${\PPPP}\rar \PP^{2\dual}$ sending each sheaf to its support.
For $t\in\PP^{2\dual}$, 
let $C_t=\phi^{-1}(t)$, $E_t=C_t/\tau$, and $\rho_t=\rho|_{C_t}:C_t\rar E_t$, where
$\phi:\CCC\rar \PP^{2\dual}$ is the natural map. We call the fiber
$P_t=(f^{2m}_{\PPPP})^{-1}(t)$ of the support map the compactified Prymian
of the pair $(C_t,\tau)$. In this section,
we will describe the structure of $P_t$ for all irreducible singular members
$C_t$ of the linear system $\CCC/\PP^{2\dual}$.

\begin{lemma}\label{degen-fibers}
Let us assume that $S$ is generic, that is, the curves $B_0\in
|\OOO_{\PP^2}(4)|$
and $\Delta_0\in |-2K_X|$ are generic. Let
$\bar{\Delta}_0=\mu(\Delta_0)\subset
\PP^2$, and let $B_0^\dual$, $\bar{\Delta}_0^\dual$ denote the
dual curves in $\PP^{2\dual}$. Let $T$ be the finite set of points which
are singularities of
the curve $B_0^\dual\cup\bar{\Delta}_0^\dual$.
Then the linear system $\CCC/\PP^{2\dual}$
contains the following singular members $C_t$:

(i) $C_t$ has a unique node $p$ if $t\in\bar{\Delta}_0^\dual\setminus T$; 
$p$ is  $\tau$-invariant,
and $\tau$ permutes the branches of $C_t$ at $p$.

(ii) $C_t$ has a unique cusp if $t$ is one of the $24$ cusps
of $\bar{\Delta}_0^\dual$.

(iii) $C_t$ has two nodes permuted by $\tau$ if $t\in B_0^\dual\setminus T$.

(iv) $C_t$ has two $\tau$-invariant nodes if $t$ is one of the $28$
nodes of $\bar{\Delta}_0^\dual$.

(v)  $C_t$ has two cusps permuted by $\tau$ if $t$ is one of the $24$ cusps
of $B_0^\dual$.

(vi) $C_t$ has $3$ nodes, only one of which is invariant under $\tau$,
if $t$ is one of the $128$ points of transversal intersection of
$B_0^\dual$ and $\bar{\Delta}_0^\dual$.

(vii) $C_t$ has one tacnode if $t$ is one of the $8$ points of
tangency of $B_0^\dual$, $\bar{\Delta}_0^\dual$.

(viii) $C_t$ is a union of two smooth conics meeting transversely
in $4$ points if $t$ is one of the $28$ nodes of $B_0^\dual$.
\end{lemma}

\begin{proof}
The proof is obvious. Remark that $B_0$, $\bar{\Delta}_0$
are totally tangent to each other; if $f_4(u_0,u_1,u_2)=0$,
$g_4(u_0,u_1,u_2)=0$ are their equations, then the pencil 
$\langle f_4,g_4\rangle$ contains the square $q^2$ of some
quadratic form $q$ in $u_0,u_1,u_2$. This follows from the fact that the
inverse image of $\bar{\Delta}_0$ in $X$ decomposes into two components,
one of which is ${\Delta}_0$. The number 28, resp. 24 is 
the number of bitangents, resp. flexes of a smooth plane quartic.
The eight tangency points of $B_0$, $\bar{\Delta}_0$ are sent by
the Gauss map to 8 tangency points of $B_0^\dual$, $\bar{\Delta}_0^\dual$.
As the degree of each of the two dual curves is 12, there remains 
$12^2-8\cdot 2=128$ points of transversal intersection, corresponding
to the non-tacnodal common tangents of $B_0$, $\bar{\Delta}_0$.
\end{proof}

\begin{lemma}\label{degen-pryms} 
Let $t\in\PP^{2\dual}$ and $P_t=(f^{2m}_{\PPPP})^{-1}(t)$. Assume that $C_t$ is irreducible.
Then the following assertions hold:

(i) The varieties $P_t$ constructed for different values
of $m$ are isomorphic to each other.

(ii)  $P_t$ has an action of the $2$-dimensional algebraic group
$G_t=\Prym(C_t,\tau)$, and the locus
$P_t^*$ of invertible sheaves in $P_t$ is a finite union of orbits of $G_t^\circ$
on which the action is free.
\end{lemma}

\begin{proof}
(i) There is a canonical isomorphism $\PPrym^{2m}(\CCC,\tau)\simeq \PPrym^{2m+4}(\CCC,\tau)$
given by $[\FFF]\mapsto [\FFF\otimes\OOO_S(H)]$. 
Thus it suffices to consider only the values $m=0$ and 1. In this case
there is no isomorphism of the relative Prymians, but there are noncanonical
isomorphisms of the individual fibers $P_t^0=(f^0_{\PPPP})^{-1}(t)\simeq
P_t^2=(f^2_{\PPPP})^{-1}(t)$. Such an isomorphism can be associated to any of the
$\tau$-invariant theta-characteristics $\theta$ of $C_t$ by 
$[\FFF]\mapsto [\FFF\otimes\theta]$. One can choose $\theta=\rho_t^{-1}(\xi)$,
where $\xi$ is a ramification point of the double cover
$\mu_t:E_t\rar\ell_t\simeq\PP^1$,
and $\mu_t=\mu|_{E_t}$.

(ii) In the case when $\LLL$ is an invertible sheaf on $C_t$ and $\FFF$ is a
\mbox{rank-1} torsion-free sheaf, we have $\tau (\FFF\otimes\LLL)=\tau (\FFF)\otimes\tau (\LLL)$,
and similarly for $\iota$. This implies that $G_t$ acts on $P_t$ by tensor multiplication
of the corresponding sheaves. The action is obviously free
on $P_t^*$.

By (i), we can assume that $m=0$.
By  \cite{AIK}, $\bar{J}(C)$ is irreducible. Further, ${\PPPP}$ is the
4-dimensional fixed locus of the involution $\kappa$ 
on $\MMM$ whose differential has exactly
2 eigenvalues $-1$ at any point of $P_t^*={\PPPP}\cap J(C_t)$, and  one of these eigenvalues
corresponds to the reflection with respect to a plane in the base $\PP^3$,
while the other to a reflection in the fiber $J(C_t)$. Thus  
every connected component of $P_t^*$ is 2-dimensional, hence
isomorphic to $G_t^\circ$.
\end{proof}

Remark that we have not proved that
$P_t$ has no 2-dimensional components contained entirely
in the non-locally-free locus. We will get this as a consequence of a case-by-case
description of a natural stratification of $P_t$ for the degenerate curves $C_t$
listed in Lemma \ref{degen-fibers}. Let us fix $t$ and omit the subscript $t$
from the symbols $C_t$, $P_t$, etc. In this section,
we consider only the case when $C$
is irreducible. By Lemma \ref{degen-pryms}, (i), $P$ does not depend
on $m$, so we can assume $m=0$.
According to \cite{Cook-1}, $\bar{J}(C)$
has a stratification in smooth strata whose codimension is equal to
the index $i(\FFF)$ of the sheaves $\FFF$
represented by the points of these strata.
The index is defined as follows. Let $\nu:\tilde C\rar C$ be the normalization map.
Then there exists a factorization $\tilde C\xrightarrow{\nu''} C'\xrightarrow{\nu'} C$ of $\nu$
such that $\nu^{\prime *}(\FFF)/\mathrm{(torsion)}$ is invertible, and
$i(\FFF)$ is the minimum of $\length (\nu'_*\OOO_{C'}/\OOO_C)$ over 
such factorizations. The index takes values between 0 and $\delta(C)=
\length (\nu_*\OOO_{\tilde C}/\OOO_C)=p_a(C)-g(C)$, and $\FFF$ is invertible
if and only if $i(\FFF)=0$.

Let $J_i(C)$ be the union of strata of
codimension $i$ in $\bar{J}(C)$ ($0\leq \mbox{$i\leq 3$}$); obviously, $J_0(C)=J(C)$. 
We will
denote by $P_i$ the intersection $J_i(C)\cap P$. 
Then $P_0=\Prym (C,\tau)$ is an algebraic
group of dimension~2, which we denoted by $G_t$ in Lemma \ref{degen-pryms}.
As we will see, for $i>0$, the value of
$i$ may differ from the codimension of $P_i$ in $P$.
We will determine the strata $P_i$ for all the singular
members of the linear system  $\CCC/\PP^{2\dual}$.

\begin{proposition}\label{prym-strata}
Assume that $m=0$.
In the notation of Lemma \ref{degen-fibers}, 
let $t\in B_0^\dual\cup\bar{\Delta}_0^\dual$. Denote by
$\nu:\tilde C\rar C$ the normalization of $C=C_t$. The map
$[\FFF]\mapsto [\nu^{*}(\FFF)/\mathrm{(torsion)}]$, when restricted
to $J_i(C)$, is a morphism $J_i(C)\rar \Pic^{-i}(\tilde C)$,
which will be denoted by $\ups_i$. The involution $\tau$ lifts to
$\tilde C$ and to some partial normalizations of $C$, and we will
use the same symbol $\tau$ to denote such a lift.

In the first seven cases of Lemma \ref{degen-fibers}, all the
nonempty strata $P_i$ are described as follows:

(i) $P=P_0\sqcup P_1$, $P_0=\ups_0^{-1}(\Prym (\tilde C, \tau))$,\ \ $P_1\simeq \Prym (\tilde C, \tau)$. Here $\Prym (\tilde C, \tau)$ is an elliptic curve lying
in the abelian surface $J(\tilde C)$, and 
$\ups_0:J(C)\rar J(\tilde C)$ is a group morphism with kernel $\CC^*$.
Thus $P_1$ is an elliptic curve, and $P_0$ is an extension
of an elliptic curve by $\CC^*$.

(ii) $P=P_0\sqcup P_1$, $P_0=\ups_0^{-1}(\Prym (\tilde C, \tau))$,\ \ $P_1\simeq \Prym (\tilde C, \tau)$ as in~(i), but now $\ker\ups_0\simeq\CC$ and $P_0$ is an extension of an elliptic curve by the additive group $\CC$.

(iii) $P=P_0\sqcup P_2$, $P_0$ is a $\CC^*$-extension of the elliptic
curve $J(\tilde C)\simeq\tilde C$, and $P_2\simeq J(\tilde C)$ (thus $\codim_PP_2=1$).

(iv) $P=P_0\sqcup P_1\sqcup P_2$, $P_0\simeq\bigsqcup\limits_{k=1}^4
\CC^*\times\CC^*$, $P_1\simeq \bigsqcup\limits_{k=1}^8 \CC^*$, and
$P_2$ is a finite set, consisting of $4$ points.

(v) $P=P_0\sqcup P_2$, $P_0$ is a $\CC$-extension of the elliptic
curve $J(\tilde C)\simeq\tilde C$, and $P_2\simeq J(\tilde C)$ (thus $\codim_PP_2=1$).

(vi) $P=P_0\sqcup \ldots\sqcup P_3$, $P_0\simeq \CC^*\times\CC^*$, $P_1\simeq 
P_2\simeq \CC^*$, and $P_3$ is one point, corresponding to the sheaf
$\nu_*(\OOO_{\PP^1}(-3))$, where we have identified $\tilde C$ with $\PP^1$.

(vii) $P=P_0\sqcup P_2$, $P_0$ is an irreducible extension of the elliptic curve $J(\tilde C)\simeq\tilde C$ by $\CC\times (\ZZ/2\ZZ)$, 
and $P_2\simeq\tilde C$.
\end{proposition}

\begin{proof}
(i) Here $\delta(C)=1$, so $\bar J(C)=J_0(C)\sqcup J_1(C)$,
hence $P=P_0\sqcup P_1$.
Let $\{p',p''\}=\nu^{-1}(p)$, $\{ p\}=\Sing C$. If $\FFF\in J_0(C)$,
then $\LLL=\nu^*(\FFF)$ is invertible and $\FFF$ is a subsheaf of
$\nu_*(\LLL)$ with quotient $\CC_p$. For fixed $\LLL$, the invertible
subsheaves of $\nu_*(\LLL)$ with quotient $\CC_p$ are parametrized
by $\CC^*$. One can easily describe them in terms of the
corresponding line bundles. We will not distinguish in the notation the invertible
sheaves and the corresponding line bundles. So, we represent
the line bundle $\FFF$ as the result of gluing together the fibers $\LLL_{,p'},
\LLL_{,p''}$ of the line bundle $\LLL$. To this end, we choose some rational
section $e$ of $\LLL$ which trivializes $\LLL$ on an open set containing
both $p'$ and $p''$, and then the gluing is determined by a factor $\lambda
\in\CC^*$ by the following rule: $\FFF=\LLL/(e_{,p'}\sim \lambda e_{,p''})$.
Let us denote it by $\LLL[_e^\lambda]$.
In the language of sheaves, $\FFF=\LLL[_e^\lambda]$ is described as the subsheaf
of $\nu_*(\LLL)$ with the following stalks:
$$
\FFF_z=(\nu_*(\LLL))_z\ \ \forall\ \ z\in C\setminus\{p\},\ \ 
\mbox{and} \ \ \FFF_p=\OOO_p\cdot (e_{p'}+ \lambda e_{p''})\ .
$$

Now we will seak the triples $(\LLL,e,\lambda)$ for which $\LLL[_e^\lambda]
\in P$, that is, $\tau^*(\LLL[_e^\lambda])\otimes\LLL[_e^\lambda]\simeq
\OOO_C$.
The involution $\tau:\tilde C\rar\tilde C$ lifts in a natural way to a map
$\tau_\#:\tau^*(\LLL)\rar \LLL$, understood as
a map of the total spaces either of sheaves, or of line bundles,
and we choose $e^\tau:=\tau_\#^{-1}(e)$ as a trivialization of
$\tau^*(\LLL)$ in the neighborhood of $p',p''$. As $\tau$ permutes
$p',p''$, we have
$\tau^*(\LLL[_e^\lambda])=\tau^*(\LLL)[_{e^\tau}^{\lambda^{-1}}]$
and $\tau^*(\LLL[_e^\lambda])\otimes\LLL[_e^\lambda]=
(\tau^*(\LLL)\otimes\LLL)\genfrac{[}{]}{0pt}{}{1}{e^\tau\otimes e}$.
Thus the necessary condition for $\LLL[_e^\lambda]
\in P$ is $\tau^*(\LLL)\otimes\LLL\simeq \OOO_{\tilde C}$.
Assume it is satisfied, then fix an isomorphism and denote by
$g$ the image of $e^\tau\otimes e$ in $\OOO_{\tilde C}$.
Then $(\tau^*(\LLL)\otimes\LLL)\genfrac{[}{]}{0pt}{}{1}{e^\tau\otimes e}
\simeq \OOO_{\tilde C}\genfrac{[}{]}{0pt}{}{g(p'')/g(p')}{1}$.
Using the canonical isomorphisms $\tau^*(\LLL)\otimes\LLL
=\LLL\otimes\tau^*(\LLL)$ and $\tau^*(\tau^*(\LLL))=\LLL$, we may
claim that $\tau_\#\circ\tau_\#=\id_\LLL$ and that $e^\tau\otimes e$
is $\tau$-invariant. Hence $g$ is also $\tau$-invariant and 
$g(p'')/g(p')=1$. Thus $(\tau^*(\LLL)\otimes\LLL)\genfrac{[}{]}{0pt}{}{1}{e^\tau\otimes e}
\simeq\OOO_{\tilde C}\genfrac{[}{]}{0pt}{}{1}{1}=\OOO_C$,
and we conclude that $\LLL[_e^\lambda]\in P$ as soon as
$\LLL\in \Prym (\tilde C, \tau)$, and this does not depend on the
choice of $e, \lambda$. We have proved the part of (i) concerning
$P_0$. 

Let now $\FFF\in J_1(C)$. Then $\FFF\simeq\nu_*(\LLL)$ for
$\LLL\in\Pic^{-1}(\tilde C)$. To express the dual
$\FFF^\dual=\EXT^1_{\OOO_S}(\FFF,\OOO_S(-C))$ in terms of
$\LLL^{-1}$, we consider $\nu$ as an embedded resolution:
let $\sigma:\tilde S\rar S$ be the blowup at $p$,
$\tilde C$ the proper transform of $C$ in $\tilde S$, and
$\nu=\sigma|_{\tilde C}$. Then, by the relative duality
for $\sigma$ (see \cite{Ha}, p. 210),
$$
\FFF^\dual\simeq \sigma_*(
\EXT^1_{\OOO_{\tilde S}}(\LLL,\OOO_{\tilde S}(-\sigma^*(C))\otimes\omega_{\tilde S/S}))
\simeq \nu_*(\LLL^{-1}(-p'-p'')).
$$
We have:
$$
\begin{array}{ll}
& [\FFF]\in P_1=\Fix (\kappa)\cap J_1(C)\\
\equi & (\tau^*\FFF)^\dual\simeq \FFF\\
\equi & (\tau^*\LLL)^{-1}(-p'-p'')\simeq \LLL\\
\equi & [\LLL (p')]\in\Prym (\tilde C, \tau).
\end{array}
$$

Thus $P_1\simeq \Prym (\tilde C, \tau)$ via the mutually
inverse maps
$$
\begin{array}{cccccccl}
P_1 &\lra& \Prym (\tilde C, \tau)& , & \Prym (\tilde C, \tau) & \lra & P_1& .
\\ 
{}[\FFF]&\longmapsto& \ups_1(\FFF)\cdot[\OOO_{\tilde C}(p')]
&& [\LLL]&\longmapsto& [\nu_*(\LLL(-p'))]&
\end{array}
$$

(ii) As in (i), $\delta (C)=1$, so $P=P_0\sqcup P_1$. Denote
by $p$ the singular point of $C$ and set $p'=\nu^{-1}(p)$.
Fix a local parameter $t$ of $\OOO_{\tilde C,p'}$ in such a way that
$\tau^*(t)=-t$, and an invertible sheaf $\LLL$ on $\tilde C$ together
with a local trivialization $e$ at $p'$. Then any invertible sheaf
$\FFF$ on $C$ such that $\nu^*(\FFF)\simeq\LLL$ can be represented as
the subsheaf of $\nu_*(\LLL)$ given by its stalks:
$$
\FFF_z=(\nu_*(\LLL))_z\ \ \forall\ \ z\in C\setminus\{p\},\ \ 
\mbox{and} \ \ \FFF_p=\OOO_p\cdot (1+bt)e_{p'}\ 
$$
for some constant $b\in\CC$. Denote this sheaf by 
$\LLL\genfrac{[}{]}{0pt}{}{b}{e;t}$. Similarly to the proof of
part (i), $\tau^*(\LLL\genfrac{[}{]}{0pt}{}{b}{e;t})=
\LLL\genfrac{[}{]}{0pt}{}{-b}{e^\tau\!;t}$, and
one easily verifies that $\LLL\genfrac{[}{]}{0pt}{}{b}{e;t}
\in P_0$ if and only if $\LLL\in \Prym (\tilde C, \tau)$ independently
of the choice of $e,t,b$. This implies the assertion about $P_0$.
The proof for $P_1$ is based on the formula for the dual
$(\nu_*(\LLL))^\dual\simeq \nu_*(\LLL^{-1}(-2p'))$, which implies
that $\nu_*(\LLL)\in P_1$\ $\Leftrightarrow$\ $\LLL(p')\in \Prym (\tilde C, \tau)$.

(iii) Here $\delta(C)=2$ and $C$ has two singular points $p_1,p_2$
permuted by~$\tau$. Let
$\nu^{-1}(p_i)=\{p_i',p_i''\}$. We can choose the notation
in such a way that $\tau (p_1')=p_2'$. $\bar J(C)$ has three strata, but
$P_1=J_1\cap P=\empt$, because if a sheaf belonging to $P$ is not
locally free at $p_1$, it is not locally  free at $p_2$ either.
Hence $P=P_0\sqcup P_2$. Let $\LLL$ be an invertible sheaf on $\tilde C$,
$e$~its rational section, regular at $p_i',p_i''$, and $(\lambda_1,
\lambda_2)\in \CC^*\times\CC^*$. Then we define $\LLL
\genfrac{[}{]}{0pt}{}{\lambda_1,\lambda_2}{e}$ as the subsheaf $\FFF$
of $\nu_*(\LLL)$ which coincides with $\nu_*(\LLL)$ over
$C\setminus \{p_1,p_2\}$ and such that
$\FFF_{p_i}=\OOO_{p_i}\cdot (e_{p'_i}+\lambda_i e_{p''_i})$, $i=1,2$.
Similarly to the part (i), one easily verifies that 
$\tau^*(\LLL
\genfrac{[}{]}{0pt}{}{\lambda_1,\lambda_2}{e})=
(\tau^*(\LLL))\genfrac{[}{]}{0pt}{}{\lambda_2,\lambda_1}{e^\tau}$.
This implies that
$\LLL
\genfrac{[}{]}{0pt}{}{\lambda_1,\lambda_2}{e}\in P_0$ if and only if
$\LLL\in \Prym (\tilde C, \tau)$ and $\lambda_1\lambda_2=1$.
Here $\tilde C$ is elliptic
and $\tilde C/\tau\simeq\PP^1$, hence $\Prym (\tilde C, \tau)
\simeq J(\tilde C)\simeq \tilde C$, and thus $P_0$
is an extension of $\tilde C$ by $\CC^*$.

Further, any sheaf from $J_2(C)$ is of the form
$\FFF=\nu_*(\LLL)$ for $\LLL\in\Pic^{-2}(\tilde C)$, and
its dual is given by $\FFF^\dual \simeq\nu_*(\LLL^{-1}(-p_1'-p_1''-p_2'-p_2''))$.
This implies that $P_2$ consists of the sheaves $\nu_*(\JJJ(-p_1'-p_1''))$,
where $\JJJ$ runs over $J(\tilde C)$.

(iv) Here $\delta(C)=2$ and $C$ has two $\tau$-invariant nodes $p_1,p_2$,
in which $\tau$ permutes the branches. Let $\nu^{-1}(p_i)=\{p_i',p_i''\}$. 
When lifted to $\tilde C$,
$\tau$ is fixed point free. As $\tau(p_i')=p_i''$, $\tau$ is a translation
on $\tilde C$ by a point of order two $[p_1''-p_1']=[p_2''-p_2']$, so that
$C/\tau=\tilde C/\tau=E$ is an elliptic curve. $\bar J(C)$ has three strata,
and we first consider a sheaf $\FFF\in J_0(C)$. As in $(iii)$, we represent
it in the form $\FFF=\LLL
\genfrac{[}{]}{0pt}{}{\lambda_1,\lambda_2}{e}$ with $(\lambda_1,
\lambda_2)\in \CC^*\times\CC^*$. Then $\tau^*(\FFF)=
(\tau^*(\LLL))\genfrac{[}{]}{0pt}{}{\lambda_1^{-1},\lambda_2^{-1}}{e^\tau}$
and $\tau^*(\FFF)\otimes\FFF=(\tau^*(\LLL)\otimes\LLL)
\genfrac{[}{]}{0pt}{}{1\ ,\ 1}{e^\tau\otimes e}$. On $\tilde C$,
$\tau^*(\LLL)\simeq\LLL$ for any degree-0 invertible $\LLL$, so
$\LLL\in P_0$ if and only if $\LLL$ is of order 2 in $J(\tilde C)$.
Thus $P_0=\ups_0^{-1}(J(\tilde C)_{(2)})$ is the disjoint union of 4 copies
of $\CC^*\times\CC^*$.

$J_1(C)$ consists of the subsheaves $\FFF\subset\nu_*(\LLL)$ with
$\LLL\in\Pic^{-1}(\tilde C)$, which coincide with $\nu_*(\LLL)$
over $C\setminus \{p_i\}$ and such that $\FFF_{p_i}=\mbox{$\OOO_{p_i}\cdot (e_{p'_i}+\lambda e_{p''_i})$}$ for one of the values of $i=1$ or 2
($\lambda\in\CC^*$).
Let us denote such a sheaf by $\LLL
\genfrac{[}{]}{0pt}{}{\lambda,\, \cdot\,}{e}$ if $i=1$
and $\LLL
\genfrac{[}{]}{0pt}{}{\, \cdot\, ,\lambda}{e}$ if $i=2$.
Let, for example, $\FFF=\LLL
\genfrac{[}{]}{0pt}{}{\lambda,\, \cdot\,}{e}$. We have $\FFF^\dual\simeq
(\LLL^{-1}(-p_2'-p_2''))
\genfrac{[}{]}{0pt}{}{\lambda^{-1},\, \cdot\,}{e\sdual}$,
so that $\kappa (\FFF)=
(\tau^*(\FFF))^\dual\simeq
\bigl(\tau^*(\LLL)^{-1}(-p_2'-p_2'')\bigr)
\genfrac{[}{]}{0pt}{}{\lambda,\, \cdot\,}{(e^\tau)\sdual}$.
A necessary condition for $\FFF\in P_1$ is $\tau^*(\LLL)^{-1}(-p_2'-p_2'')
\simeq \LLL$, or equivalently, $(\tau^*(\LLL(p_2')))^{-1}\simeq \LLL(p_2')$.
Let it be satisfied, and let us fix such an
isomorphism. Via a natural embedding $\LLL\into \LLL(p_2')$, we can consider
$e$ as a rational section of $\LLL(p_2')$, regular and nonvanishing
at $p_1',p_1''$, and the latter isomorphism
sends $(e^\tau)^\dual$ to $ge$ for some $g\in \CC(C)$. Then
$$
\begin{array}{ll}
&\FFF\in P_1 \\ \Leftrightarrow & 
(\tau^*(\LLL(p_2')))^{-1}\genfrac{[}{]}{0pt}{}{\lambda,\, \cdot\,}{(e^\tau)\sdual}
\simeq \LLL(p_2')\genfrac{[}{]}{0pt}{}{\lambda,\, \cdot\,}{e}\\
\Leftrightarrow & 
\LLL(p_2')\left[
\genfrac{}{}{0pt}{}{\frac{g(p_1'')}{g(p_1')}\textstyle \lambda\textstyle ,\ \textstyle \cdot\  }{\textstyle e}
\right]
\simeq \LLL(p_2')\genfrac{[}{]}{0pt}{}{\lambda,\, \cdot\,}{e}\\
\Leftrightarrow & g(p_1')=g(p_1'')
\end{array}
$$

It is easily seen that $g$ is $\tau$-invariant, so the last condition
is satisfied. Thus $\LLL\genfrac{[}{]}{0pt}{}{\lambda,\, \cdot\,}{e}\in P_1$
if and only if 
$\LLL(p_2')$ is one of the 4 points of second order in $J(\tilde C)$,
independently of $e,\lambda$, and this gives 4 components of $P_1$,
each isomorphic to $\CC^*$. The other four are given by the
sheaves $\LLL
\genfrac{[}{]}{0pt}{}{\, \cdot\, ,\lambda}{e}$
for which $\LLL(p_1')$ is a point of second order in $J(\tilde C)$.

Finally, $P_2$ consists of $4$ sheaves $\nu_*\LLL$, for which
$\LLL(p_1'+p_2')\in J(\tilde C)_{(2)}$.

(v) As in (iii), $P_1=\empt$. Denote $\Sing C=\{p_1,p_2\}$,
$p_i'=\nu^{-1}(p_i)$. We represent the sheaves from $J_0(C)$
in the form $\LLL\genfrac{[}{]}{0pt}{}{b_1,b_2}{e;t_1,t_2}$,
where $\LLL$ runs over $J(\tilde C)$, $e$ is a rational section
of $\LLL$ trivializing it at $p_1',p_2'$, and $t_i$ are
local parameters at $p_i'$ such that $\tau^*(t_i)=t_{3-i}$ ($i=1,2$).
The sheaf $\FFF=\LLL\genfrac{[}{]}{0pt}{}{b_1,b_2}{e;t_1,t_2}$
is defined as the subsheaf of $\nu_*\LLL$ which coincides with
$\nu_*\LLL$ over $C\setminus\{p_1,p_2\}$ and such that
$\FFF_{p_i}=\OOO_{p_i}\cdot(1+b_it_i)e_{p_i'}$ for $i=1,2$.
Then $\tau^*(\FFF)=\tau^*(\LLL)\genfrac{[}{]}{0pt}{}{b_2,b_1}{e^\tau;t_1,t_2}$,
and the stalk of $\tau^*(\FFF)\otimes\FFF$ at $p_i$, as an $\OOO_{p_i}$-submodule
of the stalk of $\nu_*(\tau^*(\LLL)\otimes\LLL)$, is generated by
$(1+b_1t_i)(1+b_2t_i)e_{p_i'}^{\tau}\otimes e_{p_i'}$. As $t_i^2\in\mathfrak m_{p_i}$, we conclude that
$\tau^*(\FFF)\otimes\FFF=(\tau^*(\LLL)\otimes\LLL)
\genfrac{[}{]}{0pt}{}{b_1+b_2,b_1+b_2}{e^\tau\otimes e;t_1,t_2}$
and that $\FFF\in P_0$ \ $\Leftrightarrow$ \ $b_1+b_2=0$. Thus
$P_0$ is a $\CC$-extension of $J(\tilde C)$.

The stratum $J_2(C)$ consists of the sheaves $\nu_*(\LLL)$, where
$\LLL$ runs over $\Pic^{-2}(\tilde C)$, and
$(\nu_*(\LLL))^\dual \simeq\nu_*(\LLL^{-1}(-2p_1'-2p_2'))$.
This implies that $P_2$ consists of the sheaves $\nu_*(\JJJ(-p_1'-p_2'))$,
where $\JJJ$ runs over $J(\tilde C)$, hence $P_2\simeq \tilde C$.

(vi)  Here $\delta(C)=3$ and $\bar J(C)$ has 4 strata.
We set $\Sing C=\{p_0,p_1,p_2\}$, $\tau (p_0)=p_0$, $\tau(p_1)=p_2$,
$\nu^{-1}(p_i)=\{p_i',p_i''\}$, $\tau (p_0')=p_0''$, $\tau(p_1')=p_2'$.
As $\tilde C\simeq\PP^1$, the open stratum $J_0(C)$ consists of the subsheaves
$\OOO_{\tilde C}\genfrac{[}{]}{0pt}{}{
\lambda_0,\lambda_1,\lambda_2}{1}$ of $\nu_*\OOO_{\tilde C}$
which coincide with $\nu_*\OOO_{\tilde C}$ over $C_\ns$ and whose
stalk at $p_i$ is generated by $1_{p_i'}+\lambda_i 1_{p_i''}$.
We have $\tau^*\FFF\simeq \OOO_{\tilde C}\genfrac{[}{]}{0pt}{}{
\lambda_0^{-1},\lambda_2,\lambda_1}{1}$ and
$\tau^*\FFF\otimes\FFF\simeq \OOO_{\tilde C}\genfrac{[}{]}{0pt}{}{
1,\lambda_1\lambda_2,\lambda_1\lambda_2}{1}$. Thus $P_0\simeq (\CC^*)^2$
is the subtorus of $J(C)\simeq (\CC^*)^3$ singled out by the equation
$\lambda_1\lambda_2=1$. Similarly, we can describe the other strata:
$$
\begin{array}{l}
P_1=\left\{\OOO_{\tilde C}(-pt)\genfrac{[}{]}{0pt}{}{
\cdot,\lambda,\lambda^{-1}}{1}\right\}_{\lambda\in\CC^*}\simeq\CC^*;\\
P_2=\left\{\OOO_{\tilde C}(-2pt)\genfrac{[}{]}{0pt}{}{
\lambda,\cdot,\cdot}{1}\right\}_{\lambda\in\CC^*}\simeq\CC^*;\\
P_3=\{\nu_*(\OOO_{\tilde C}(-3pt))\}\ = \ \mbox{1 point.}
\end{array}
$$
To clarify the notation, we remind that $\tilde C\simeq\PP^1$,
so that $\tau$ has two fixed points
on $\tilde C$. We use one of them, denoted $pt$, to embed the sheaves
$\nu^*\FFF/\mbox{(torsion)}$ with $\FFF\in P_i$ into $\OOO_{\tilde C}$ as the $\tau$-invariant subsheaves
$\OOO_{\tilde C}(-i\:pt)$, and the section 1 of $\OOO_C$ is considered
as a rational trivialization of $\OOO_{\tilde C}(-i\:pt)$.

(vii) Here $C$ has one tacnode $p$, $\delta (C)=2$,
so that $\bar J(C)$ has three strata. Denote by $p_1,p_2$ the preimages
of $p$ in $\tilde C$ and fix some local parameters $t_i$ at $p_i$.
The points $p_i$ are $\tau$-invariant, and we can choose the $t_i$
in such a way that $\tau^*(t_i)=-t_i$. We will identify the formal completion
$B=(\nu_*\OOO_{\tilde C})_p\!\!\hat{}$ \ \ of $\nu_*\OOO_{\tilde C}$ at $p$ with
$\CC[[t_1]]\times\CC[[t_2]]$. We can further restrict the
choice of the $t_i$ so that the formal completion $A=\OOO_p\hat{}$ is given
by
$$A=\{(a_0+a_1t_1+a_2t_1^2+\dots,
b_0+b_1t_2+b_2t_2^2+\dots)\in B\mid a_0=b_0,\ a_1=b_1\}.
$$
Denote by $\mathfrak c$ the conductor of $\OOO_p$ in $(\nu_*\OOO_{\tilde C})_p$:
$$
\mathfrak c=\{u\in \OOO_p\mid u(\nu_*\OOO_{\tilde C})_p\subset \OOO_p\} .
$$
For its completion, we have $\hat{\mathfrak c}=A(t_1^2,0)+A(0,t_2^2)=
B(t_1^2,t_2^2)$.

The description of $\bar J(C)$ that we will expose here is
similar to that given in \cite{Cook-1}.
To each $\FFF\in\bar J(C)$,
we have assigned the invertible sheaf $\LLL=\nu^*\FFF/\mbox{(torsion)}$
on $\tilde C$. Let us fix a local trivialization of $\LLL$
by a rational section $e$, regular and nonvanishing
at $p_1,p_2$. Then $\FFF$ can be described as a subsheaf of
$\nu_*\LLL$ which coincides with $\nu_*\LLL$ out of $p$ and such that
$\FFF_p\subset (\nu_*\LLL)_p$ is an $\OOO_p$-submodule of colength $2-i(\FFF)$.

Consider the case when $\FFF\in J_0(C)$. Then $\LLL$ is of degree $0$
and $\FFF_p$ is of colength 2. Quotienting by $\mathfrak c$,
we obtain the vector plane $\FFF_p/\mathfrak c\FFF_p$ 
in the 4-dimensional vector space
$V=(\nu_*\LLL)_p/\mathfrak c(\nu_*\LLL)_p$.
This gives a one-to-one correspondence between the 
sheaves $\FFF\in J_0(C)$ with the same assigned $\LLL$
and the vector planes in $V$ which are principal $\OOO_p/\mathfrak c$-modules.
Such planes form a locally closed subset $U_\LLL$ of the Grassmannian $G(2,V)$,
and to describe it, we can go over to the formal completions.
Let $\bar A=A/\hat{\mathfrak c}$, $\bar B=B/\hat{\mathfrak c}$. Using
$e$ as a generator of $(\nu_*\LLL)_p$, we will identify
$V$ with $\bar B$.

Thus we can choose $(1,0),(\bar t_1,0),(0,1),(0,\bar t_2)$ as a basis of $V$,
where the bar over $t_i$ means taking the coset modulo $\hat{\mathfrak c}$.
Then the \mbox{2-planes} in $\bar B$, invariant under the multiplication
by the elements of $\bar A=\langle (1,1),(\bar t_1,\bar t_2)\rangle$,
form a 2-dimensional quadratic cone $Q_\LLL$. If we introduce the
Pl\"ucker coordinates $p_{ij}$ associated to the above basis
of $\bar B$, then $G(2,\bar B)=G(2,4)$ is the 
Pl\"ucker quadric in $\PP^5$ with equation 
$p_{12}p_{34}-p_{13}p_{24}+p_{14}p_{23}=0$, and $Q_\LLL$ is
the linear section of $G(2,4)$ defined by $p_{13}=p_{14}+p_{23}=0$.
The 2-planes that are principal \mbox{$\bar A$-modules} are parametrized
by the open subset of $Q_\LLL$, the complement of two generators
of the cone: $U_\LLL=Q_\LLL\setminus (\ell_1\cup\ell_2)$. If we denote
by $e_{ij}$ the point of $G(2,4)$ for which $p_{ij}=1$ and all the other
$p_{kl}$ are zero, then $\ell_1=\langle e_{12}, e_{24}\rangle$
and $\ell_2=\langle e_{24}, e_{34}\rangle$. The following
map is an isomorphic parametrization of $U_\LLL$:
$$
\Pi:\CC^*\times\CC\lra U_\LLL, \ \ \ \ (\lambda,b)\mapsto
[\bar A\cdot (1,\lambda +b\bar t_2)],
$$
or in Pl\"ucker coordinates,
$$
(p_{12}:p_{13}:p_{14}:p_{23}:p_{24}:p_{34})=(1:0:\lambda:-\lambda:-b:\lambda^2).
$$
Remark, that $U_\LLL$ is an orbit of the group of units $\bar B^\times$ and
does not depend on the choice of $e$, for different $e$'s differ
by a unit of $(\nu_*\OOO_{\tilde C})_p$. Thus $J_0(C)$ is a $\CC^*\times\CC$-bundle
over $J(\tilde C)$. We will denote the invertible sheaf on $C$
corresponding to the plane $\Pi(\lambda,b)$ by $\LLL
\genfrac{[}{]}{0pt}{}{\lambda,b}{e;t_2}$.

Now we compute the acton of $\kappa=\iota\circ\tau^*$ on $U_\LLL$:
$$
\LLL\genfrac{[}{]}{0pt}{1}{\lambda,b}{e;t_2}\sendsto{\tau^*}
(\tau^*\LLL)\genfrac{[}{]}{0pt}{1}{\lambda,-b}{e^\tau;t_2}\sendsto{\iota}
(\tau^*\LLL)^{-1}\genfrac{[}{]}{0pt}{1}{\lambda^{-1},\, b\lambda^{-2}}{(e^\tau)\sdual\, ;\ t_2}\ .
$$
As $\tilde C/\tau\simeq\PP^1$, we have $(\tau^*\LLL)^{-1}\simeq\LLL$ for
any $\LLL\in J(\tilde C)$. Fix such an isomorphism; then it sends 
$(e^\tau)\sdual$ to $ge$ for some $g\in\CC(C)$. The $\tau$-invariance
of $g$ implies that it has no linear terms in $t_i$, and we deduce that
$\kappa \left(\LLL\genfrac{[}{]}{0pt}{1}{\lambda,b}{e;t_2}\right)=
\LLL\genfrac{[}{]}{0pt}{1}{a\lambda^{-1},\, ab\lambda^{-2}}{e\ \ ;\ \ t_2}$,
where $a=\frac{g(p_1)}{g(p_2)}\neq 0$. Hence the fixed locus of
$\kappa$ in $U_\LLL$ is given by $\lambda=\pm\sqrt{a}$, which singles out
two generators of the cone (with deleted vertex). Remark, that
$\kappa$ is a restriction of a linear involution on $\PP^5$:
\begin{multline}\label{kappa-linear}
(p_{12}:p_{13}:p_{14}:p_{23}:p_{24}:p_{34})\sendsto{\kappa}\\
(p_{34}:p_{13}:ap_{14}:ap_{23}:ap_{24}:a^2p_{12})
\end{multline}

Thus $P_0$ is fibered over $J(\tilde C)$ with each fiber the disjoint
union of two copies of $\CC$. To see that the family of components
of the fibers is an irreducible double cover of $J(\tilde C)$, one can
argue as follows. Write down the double cover $\tilde C\rar\tilde C/\tau
\simeq\PP^1$ in coordinates:
$$
\tilde C=\{y^2=(x-x_1)(x-x_2)(x-x_3)\} \ \rar \PP^1, \ \ 
(x,y) \mapsto\ x.
$$
It is ramified at the 4 points $p_i=(x_i,0)$ ($i=1,2,3$) and $p_4=\infty$. Parametrize $J(\tilde C)$ by the map $\tilde C\rar J(\tilde C)$,\ \  $q\mapsto
[\OOO(q-p_4)]$. Embed $\OOO(q-p_4)$ into the constant sheaf $\CC(C)$
in the natural way and use $1\in\CC(C)$ as the rational trivialization
$e$ of $\LLL$. Then the function $g$ introduced in the previous paragraph
is given by $g=x-x(q)$, where $q=(x(q),y(q))$, and
the equation $\lambda^2=\frac{g(p_1)}{g(p_2)}$, whose two solutions
provide two components of the fiber of the fibration $f: P_0\rar J(\tilde C)$ over $[\OOO(q-p_4)]$, becomes
\mbox{$\lambda^2=\frac{x_1-x(q)}{x_2-x(q)}$.} Varying $x=x(q)$, we obtain
the curve $\Gamma$ with equation \mbox{$\lambda^2=\frac{x_1-x}{x_2-x}$},
and the connected components of fibers of $f$ are
parametrized by the normalization of $\Gamma\times_{\PP^1}
\tilde C$. The latter is a nonramified double cover of $\tilde C$.

Now we will determine the lower-dimensional strata of $P$.
Instead of looking for the non-invertible sheaves $\FFF$ in $\bar J(C)$
as $\OOO_C$-submodules of colength $2-i$ in $\nu_*\LLL$ with
$\deg\LLL=-i$, we can get all of them as $\OOO_C$-submodules of colength $2$ with $\LLL$ of degree 0, parametrized by the points of $Q_\LLL\setminus U_\LLL$.
In fact, it is easy to see that the cones $Q_\LLL$ fit into
an algebraic family over $J(\tilde C)$ and that this family is the normalization
of $\bar J(C)$ (see \cite{Cook-2}), thus any non-invertible sheaf in $\bar J(C)$
is in the closure of $U_\LLL$ for some $\LLL\in J(\tilde C)$.

As follows from (\ref{kappa-linear}), $\kappa$ permutes $\ell_1$ and $\ell_2$,
thus the only fixed point of $\kappa$ in $Q_\LLL\setminus U_\LLL$
is the vertex $e_{24}$ of the cone. It corresponds to the 2-plane
$\langle (\bar t_1,0), (0,\bar t_2)\rangle$ in $\bar B$. Thus the
associated sheaf $\FFF$ has for its stalk at $p$
$$
\FFF_p=\OOO_p\cdot t_1e_{p_1}+\OOO_p\cdot t_2e_{p_2}=
(\nu_*(\LLL(-p_1-p_2)))_p,
$$
and as $\FFF$, $\nu_*\LLL$, $\nu_*(\LLL(-p_1-p_2))$ coincide
on $C\setminus\{p\}$, we conclude that $\FFF\simeq
\nu_*(\LLL(-p_1-p_2))$. It is of index 2, and we see that
$P_2\simeq J(\tilde C)$, $P_1=\empt $. This ends the proof of the proposition.
\end{proof}

\section{Further properties of ${\PPPP}^{2m}$}

Fujiki has constructed a number of irreducible symplectic \mbox{$V$-manifolds}
of dimension 4 with at worst isolated singularities
as partial desingularizations of a finite
quotient of the product of two symplectic surfaces.
Among his examples, there are two
with 28 singular points of the same type that the singular points of
${\PPPP^{2m}}$, see Table 1 on p. 225 and Remark 13.2.4 on p. 227 of
\cite{F}. 

These two examples are obtained by the following construction.
Let $H$ be a finite group of symplectic automorphisms of a K3 surface $S$,
and $\theta\in \Aut H$ such that $\theta^2=\id$. Then $H$ acts
on $S\times S$ by the rule $h:(s,t)\mapsto (hs,\theta(h)t)$.
Define $G\subset \Aut(S\times S)$ as the subgroup generated by
$H$ and the involution $(s,t)\mapsto (t,s)$. Then $K=S\times S/G$
is a symplectic $V$-manifold, in general, with non-isolated singularities.
The two examples under consideration correspond to $H=\ZZ/2\ZZ$
or $H=(\ZZ/2\ZZ)^3$ and $\theta:h\mapsto h^{-1}$. For these
$H,\theta$, the blowup of the 2-dimensional components of the singular
locus of $K$ provides two irreducible symplectic $V$-manifolds $Y_1$, $Y_2$
with 28 singular points of analytic type of
$(\CC^4/\{\pm 1\},0)$. They have the same 
Euler characteristic and the Hodge numbers. The symmetries for
the Hodge diamond of a symplectic $V$-manifold imply that the whole
Hodge diamond of $Y_i$ is determined by the three of them,
$h^{1,1}=14$, $h^{1,2}=0$, $h^{2,2}=162$, and the Euler number
is $\chi (Y_i)=8+4h^{1,1}+h^{2,2}-4h^{1,2}=226$.

The easiest way to prove that Fujiki's examples are
different from ${\PPPP}^{2m}$ is to compute the Euler number. 
Recall that there are at most two non-isomorphic varieties among the
${\PPPP}^{2m}$: \ ${\PPPP}^{0}$ and ${\PPPP}^{2}$. 

\begin{proposition} The varieties ${\PPPP}^{0}$ and ${\PPPP}^{2}$
have the same topological Euler number, equal to $268$.
\end{proposition}

\begin{proof}
Let $\PPPP$ denote either one of the varieties ${\PPPP}^{2m}$,
$f:\PPPP\rar\PP^{2\dual}$ the natural map. We introduce a stratification
$(\Pi_i)_{i=0,\ldots,8}$ of $\PP^{2\dual}$ as follows:
$\Pi_0=\PP^{2\dual}\setminus(B_0\cap\bar\Delta_0)$, the complement of
the discriminant divisor of $f$, and $\Pi_k$ for $k=1,\ldots,8$ is the locus of
points $t\in B_0\cap\bar\Delta_0$ for which the \mbox{$k$-th case} of
Lemma \ref{degen-fibers} is realized. Then we can compute
the topological Euler number of $\PPPP$ by the formula
$\chi(\PPPP)=\sum_{k=0}^8 \chi(\Pi_k)\chi(P_{t_k})$.
From Lemma \ref{prym-cc} and Proposition \ref{prym-strata},
we see that $\chi(P_{t_k})$ is the number of 0-dimensional strata
in $P_{t_k}$ and that it is different from zero only for $k=4,6,8$.
For these values of $k$, $\Pi_k$ is finite and $\chi(\Pi_k)=\#\Pi_k$.
Thus
$$
\chi(\PPPP)=28\cdot 4+128\cdot 1+28\cdot 1=268.
$$
\end{proof}

To show that ${\PPPP}^{2m}$ are irreducible symplectic $V$-manifolds
in the sense of Definition \ref{ISVM}, it remains to prove their simple
connectedness. We will start by the case $m=0$, in which we will use
a certain rational map $\Phi :S^{[2]}\dasharrow \PPPP^0$, an analog of
the Abel-Jacobi map for Prym varieties.
Recall some notation from Section \ref{section1}: $\tau:S\rar S$ is
the Galois involution of the double cover $\rho:S\rar X$,
$\mu:X\rar\PP^2$ is the double cover map, $B\subset X$ (resp. $\Delta
\subset S$) the ramification curve of $\mu$ (resp. $\rho$), $B_0=\mu(B)$,
$\Delta_0=\rho(\Delta)$. 
Let $\xi\in S^{[2]}$
be generic. Then $\xi$ is a pair of distinct points, $\xi=\{p_1,p_2\}$,
and the line $\ell_\xi =\langle\mu\rho(\xi)\rangle$ spanned by 
$\mu\rho(p_1),\mu\rho(p_2)$ in $\PP^2$ is well-defined. Let
$C_\xi=(\mu\rho)^{-1}(\ell_\xi)$. Then $\Prym (C_\xi,\tau|_{C_\xi})$
is a subvariety of $\PPPP^0$, the fiber $f^{-1}(\{\ell_\xi\})$,
where $f:\PPPP^0\rar\PP^{2\dual}$ is the natural map and $\{\ell\}$
denotes the point of $\PP^{2\dual}$ representing a line $\ell\subset\PP^2$.
Define
$$
\Phi :S^{[2]}\dasharrow \PPPP^0, \ \ \xi=\{p_1,p_2\}\mapsto
\sum_{i=1}^2[p_i-\tau(p_i)]\in \Prym (C_\xi,\tau|_{C_\xi}).
$$

Obviously, $\Phi$ is dominant. To describe the fibers of $\Phi$,
we will introduce the involution
$$
\iota_0:S^{[2]}\lra S^{[2]},\ \ \ \xi\mapsto \xi'=(\langle\xi\rangle\cap S)-\xi.
$$
Here $S$ is considered in its embedding as a quartic surface in $\PP^3$,
given by the linear system $|H|$, $\langle\xi\rangle$ stands for the
line in $\PP^3$ spanned by $\xi$, and $\xi'$ is the residual
intersection of $\langle\xi\rangle$ with $S$. This involution is regular
whenever $S$ contains no lines, which is the case for sufficiently
generic $S$ (see Lemma \ref{pairs-of-conics}). Further, $\tau$ induces on
$S^{[2]}$ an involution which
we will denote by the same symbol. 
As $\tau$ on $S$ is the restriction of a
linear involution on $\PP^3$, $\iota_0$ commutes with $\tau$, and
the composition $\iota_1=\iota_0\circ\tau$ is also an involution.

\begin{lemma}\label{fibers-of-Phi}
$\Phi$ is a rational double covering with Galois involution $\iota_1$, so that
the quotient $M=S^{[2]}/\iota_1$ is birational to $\PPPP^0$.
\end{lemma}

\begin{proof}
Let $\xi=\{p_1,p_2\}$ be generic. We have to determine all the divisors
$p_1'+p_2'$ on $C_\xi$ such that
$p_1-\tau(p_1)+p_2-\tau(p_2)\sim p_1'-\tau(p_1')+p_2'-\tau(p_2')$.
Assume this relation satisfied, and set $\delta = p_1+p_2+\tau(p_1')+\tau(p_2')$,
$\delta' =p_1'+p_2'+\tau(p_1)+\tau(p_2)$.
Then either $\delta'\neq\delta$ and $\dim |\delta|>0$,
or $\delta=\delta'$.

Let us consider the first case. There are three subcases:

(1)~$\dim |\delta|=2$. Then $\delta,\delta'\sim K=K_{C_\xi}$ are intersections
of $C_\xi$, considered as a plane quartic, with two different lines $L_1,L_2$, and
\mbox{$\tau(p_1'+p_2')$} is uniquely determined as the residual intersection
$(L_1\cap C_\xi)-p_1-p_2$. Thus there is a unique solution
$p_1'+p_2'$, different from $p_1+p_2$: $p_1'+p_2'=\tau((L_1\cap C_\xi)-p_1-p_2)=
\iota_1 (p_1+p_2)$.

(2)~$\dim |\delta|=1$ and $|\delta|$ is base point free. Then there exist
4 points $\tilde\delta$ on $C_\xi$, such that no three of them are aligned,
and $|\delta|$
consists of the residual intersections $(q\cap C_\xi)-\tilde\delta$,
where $q$ runs over the pencil of conics $|2H-\tilde\delta|$
in the plane $\langle C_\xi
\rangle$ spanned by $C_\xi$.
Remark that $\tau$ acts as a linear involution on this plane with fixed line $L_\tau$, and when $q$ runs over $|2H-\tilde\delta|$, the symmetric conic
$\tau(q)$ runs over another pencil of the same type,
$|2H-\tau(\tilde\delta)|$. As $\delta'=\tau (\delta)$,\ \  $\delta'$ belongs
to both pencils, hence the two pencils coincide. We conclude that $\tilde\delta$
is $\tau$-invariant, and hence every conic in $|2H-\tilde\delta|$ is
$\tau$-invariant. Hence $p_1'+p_2'=p_1+p_2$, which is absurd.

(3)~$\dim |\delta|=1$ and $|\delta|$ has a base point. There are two
points \mbox{$u,v\in C_\xi$} such that $|\delta|=
\{u-v+L\cap C_\xi\}$, where the line $L$ runs over the pencil $|H-v|$.
As $\delta'=\tau (\delta)$,\ \ $u, v\in L_\tau\cap C_\xi$, hence
either $\{p_1,p_2\}\cap L_\tau\cap C_\xi\neq\empt$, or
$p_1,p_2$ are aligned with one of the 4 points of $L_\tau\cap C_\xi$.
In both cases $\xi$ is non-generic, which contradicts our assumption.
\smallskip

It remains to consider the second case $\delta=\delta'$. Then $\delta$
is $\tau$-invariant, and modulo the transpositions $p_1\leftrightarrow p_2$,
$p_1'\leftrightarrow p_2'$, there are only two possibilities for which
$p_1'+p_2'\neq p_1+p_2$:

(a) $p_i'=\tau(p_i)$, $i=1,2$. Then $2(p_1+p_2)\sim 2(\tau (p_1)+
\tau(p_2))$, hence $p_1+p_2$ is nongeneric.

(b) $p_1'=p_1$, $p_2'=\tau(p_2)$. Then $2p_2\sim 2\tau(p_2)$, hence
$p_1+p_2$ is nongeneric.\smallskip

We conclude that the generic fiber of $\Phi$ consists of two elements:
$\{\xi, \iota_1(\xi)\}$.
\end{proof}

\begin{lemma}
The fixed locus $\Fix(\iota_1)$ of $\iota_1$
is the union of a nonsingular irreducible
surface $\Sigma\subset S^{[2]}$ and of $28$ isolated points.
\end{lemma}

\begin{proof}
It is obvious that the fixed point set of any biregular involution
on a smooth variety is also smooth. Consider $S$ as a quartic in
$\PP^3$. As $\tau$ has invariant curves in the linear system of
hyperplane sections $H$, it acts linearly on
$\PP^3$, and its fixed locus is the union of a plane $H_\tau$ and a
point $\infty_\tau\not\in S$. If $\iota_1(\xi)=\xi$, then
the line $\langle\xi\rangle$ is $\tau$-invariant. Hence either
$\langle\xi\rangle\subset H_\tau$, or $\langle\xi\rangle$ passes
through $\infty_\tau$.
The first case provides the 28 \ isolated
points of $\Fix(\iota_1)$, each of them being the pair of tangency
points of a bitangent to the plane quartic $\Delta_0=H_\tau\cap S$.
The second case provides the remaining part of $\Fix(\iota_1)$:
$$
\Sigma =\{\ \xi\in S^{[2]}\mid \infty_\tau\in\langle\xi\rangle,
\ \tau (\xi)\neq \xi\}.
$$

Let us call the lines through $\infty_\tau$ vertical.
A generic vertical line $L$ meets $S$ in 4 points which represent
one fiber of $\mu\rho$.
These 4 points form 6 pairs. When we vary $L$, the 6 pairs
sweep a surface $\tilde\Sigma\subset S^{[2]}$, a 6-sheeted covering of $\PP^2$.
Two of the 6 pairs are 
$\tau$-invariant, so $\tilde\Sigma$ contains an irreducible component $\Sigma_0$
which is a double covering of  $\PP^2$ and
is identified with $X=S/\tau$. 
The other 4 pairs sweep $\Sigma$, \mbox{a 4-sheeted} covering of $\PP^2$, and
we have $\tilde\Sigma=\Sigma\cup \Sigma_0$. If we assume that $\Sigma$
is reducible, then the two components of $\Sigma$ would meet along the
curve (identified with $\rho^{-1}(B)$)
of pairs of tangency points of the vertical bitangents to $S$.
This would
contradict the smoothness of $\Sigma$. Hence $\Sigma$ is irreducible.
\end{proof}

\begin{proposition}
The varieties $\PPPP^{0}$ and $M=S^{[2]}/\iota_1$ are simply connected.
\end{proposition}

\begin{proof}
We will first prove that $M$ and its resolution
of singularities $\tilde M$ are simply connected. Denote by
$\Psi$ the quotient map $S^{[2]}\rar M$.
Let $F=\Fix(\iota_1)$
and $\bar F$ the image of $F$ in $M$. Choose any point $z_0\in F$
and denote by $\bar z_0$ its image under $\Psi$.
Then any loop
based at $\bar z_0$ lifts to a loop based at $z_0$, just because $z_0$
is the unique preimage of $\bar z_0$. Hence the map
$\Psi_*:\pi_1(S^{[2]}, z_0)\rar \pi_1(M,\bar z_0)$ is surjective.
But $\pi_1(S^{[2]}, z_0)=1$, so $M$ is simply connected.

The singularities of $M$ are analytically equivalent to
$(\CC^4/\{\pm 1\}, 0)$ at the 28 isolated points of $\bar F$
and to $((\CC^2/\{\pm 1\})\times\CC^2, 0)$ along $\bar\Sigma=\Psi(\Sigma)$.
Thus a resolution of singularities can be obtained by a single
blowup $\sigma:\tilde M\rar M$ with center $\bar F$, and
the fibers of $\sigma$ over the points of $\bar F$
are the projective spaces $\PP^3$ and $\PP^1$.
Hence $\sigma$ does not change the fundamental group and $\tilde M$
is simply connected. Similarly, the blowup $\tilde\PPPP^0\rar
\PPPP^0$ of the 28 singular points
of $\PPPP^0$ pastes in 28 copies of $\PP^3$ and hence does not
change the fundamental group. We have obtained two complete smooth
varieites $\tilde\PPPP^0$, $\tilde M$ which are birational.
It follows that their fundamental groups are isomorphic. This can be deduced
from the Weak Factorization Theorem \cite{AKMW}, saying that a birational
map between complete smooth
varieites over an algebraically closed field of characteristic 0
decomposes into blowups with smoth centers or their inverses,
and from an obvious observation that a blowup of a smooth variety
with smooth center does not change the fundamental group.
We have thus proved the simple connectedness of $\PPPP^0$.
\end{proof}

\begin{lemma}\label{torsor}
Let $\GGB\subset \PPPP^0$ be the open subscheme parametrizing
invertible sheaves on the curves $C_t$, $t\in T$, where $T=\PP^{2\dual}$; it is
a group scheme over $T$ with a regular action on $\PPPP^0$.
Let $\GGGG$ denote the sheaf of cross-sections of $\GGB$
in the \'etale topology over $T$, and $\GGGG_2$ the constructible
subsheaf of $2$-torsion points. Then there exists a $1$-cocycle
$\beta$ representing an element of $H^1_{\mathrm{\acute et}}(T,\GGGG_2)$
such that $\PPPP^2\simeq \PPPP^0\times_\GGB\GGB^\beta$, where $\GGB^\beta$
is the $\GGB$-torsor defined by $\beta$.
\end{lemma}

\begin{proof}
The theta-characteristics of the curves $C_t$, that is, invertible
sheaves $\theta$ on $C_t$ such that $\theta^{\otimes 2}\simeq \omega_{C_t}$,
form a constructible sheaf $\varTheta$ with finite stalks over $T$. Let
$\varTheta^\tau$ denote the subsheaf of $\tau$-invariant theta-characteristics.
As we saw in the proofs of Lemmas \ref{prym-cc} and \ref{degen-pryms} (i),
$\varTheta^\tau$ has nonempty stalks at all the points $t\in T$.
Thus there exists an \'etale covering $(i_j:U_j\rar T)_{j\in J}$
together with local sections $\theta_j\in\Gamma(U_j,i_j^*\varTheta^\tau)$.
The translation by $\theta_j$ defines an isomorphism
$T(\theta_j):\PPPP^0_{U_j}\isoto \PPPP^2_{U_j}$. We can define the cocycle
$\beta= (\beta_{jk})$ over $U_{jk}=U_j\times_TU_k$ by
$\beta_{jk}=\pr_j^*\theta_j\otimes(\pr_k^*\theta_k)^{-1}$, where
$U_{jk}\xymatrix@1{\ar@<.6ex>[r]^{\pr_j}\ar@<-.6ex>[r]_{\pr_k}&}
\genfrac{}{}{0pt}{}{U_j}{U_k}$ are natural projections.
\end{proof}

\begin{proposition}
$\PPPP^2$ is simply connected.
\end{proposition}

\begin{proof}
Let $\PPPP$ denote either one of the varieties ${\PPPP}^{0}$ or ${\PPPP}^{2}$,
$f:\PPPP\rar\PP^{2\dual}$ the natural map, $D=B_0\cap\bar\Delta_0$
the discriminant divisor of $f$, $U=\PP^{2\dual}\setminus D$,
$E=f^{-1}(D)$, $V=\PPPP\setminus E$, so that $f_V=f|_V:V\rar U$
is a smooth projective morphism. Then $f_U$ is a locally trivial
fiber bundle in the $C^\infty$-category with general fiber $P_t$,
and there is an exact sequence of homotopy groups:
$$
\ldots\rar\pi_2(U)\xrightarrow{\partial}\pi_1(P_t)\xrightarrow{\eps}\pi_1(V)\rar\pi_1(U)\rar1 .
$$ 
It allows us to identify $\pi_1(P_t)/\im\partial$ with a subgroup
of $\pi_1(V)$. Let $(c_j)_{j\in J}$ be any generating system
for $\pi_1(P_t)/\im\partial$. Let us also fix one lift
$\tilde\gamma$ in $\pi_1(V)$ for each element $\gamma$ of
$\pi_1(U)$ different from 1. Then,
by Proposition 0.2 of \cite{Lei} and taking
into account the fact that $\pi_1(\PP^{2\dual})=1$, we obtain
a surjection $\pi_1(P_t)/\im\partial\twoheadrightarrow \pi_1(\PPPP)$
whose kernel is generated, as a normal subgroup, by the 
following set of commutators:
\begin{equation}\label{R}
R=\bigl\{ \ [\tilde\gamma, c_j]\ \mid\ j\in J,\ \gamma\in \pi_1(U)\setminus\{1\}\bigr\} .
\end{equation}

The description of $R$ in our case simplifies drastically due to the fact
that $\pi_1(P_t)\simeq\ZZ^4$ is abelian. As $\pi_1(P_t)=H_1(P_t,\ZZ)$,
the monodromy action
$M:\pi_1(U)\rar \Aut \pi_1(P_t)$ is well-defined, and for any $c\in
\pi_1(P_t)$, we have $\tilde\gamma \eps(c)\tilde\gamma^{-1}=M_\gamma(c)$
(as above, $f_{V*}(\tilde\gamma)=\gamma\in \pi_1(U)\setminus\{1\}$).
Thus we can write $\pi_1(\PPPP)\simeq \pi_1(P_t)/N$, where
$N=<\!\!R_1,R_2\!\!>_\norm$ is the normal subgroup of $\pi_1(P_t)$ generated by the
two sets of elements:\smallskip\\
\begin{minipage}{\textwidth}
\begin{itemize}
\item[\(R_1:\)] the elements of the form $M_\gamma (c_j)c_j^{-1}$, where $\gamma$
runs over \mbox{$\pi_1(U)\setminus\{1\}$}, and $(c_j)$ is a
basis of $\pi_1(P_t)$ ($j=1,\ldots,4$);
\item[\(R_2:\)] the image of any generating subset of $\pi_2(U)$.
\end{itemize}
\end{minipage}\smallskip

We will show that if $\PPPP=\PPPP^2$, then $R_1$
generates the whole of $\pi_1(P_t)$,
and thus $\pi_1(\PPPP^2)=1$. By Lemma \ref{torsor}, the smooth
locus $V=\PPPP^2_U$ of $\PPPP^2/\PP^{2\dual}$ can be obtained
by gluing together pieces $\PPPP^0_{U_j}$ of $\PPPP^0/\PP^{2\dual}$
over $U_i\cap U_j$ via transition maps which are translations
in the fibers. A~translation in a fiber induces a canonical
isomorphism of the homology groups, hence the local systems
of the groups $H_1(P_t,\ZZ)$ for $\PPPP^2_U$ and $\PPPP^0_U$
are isomorphic. Thus it suffices to see that $R_1$ generates
the whole of $H_1(P_t,\ZZ)=\pi_1(P_t)$ in the case when
$\PPPP=\PPPP^0$. Here we can use Propositioon 0.3 of op. cit.
The latter applies to the situation when $f$ has a global cross-section
meeting all the components of $E$, which is the case for 
the cross-section of neutral elements of the group scheme $\GGB$
inside $\PPPP^0$. It permits to replace the
description of the relations in the fundamental group given in (\ref{R})
by the following one:
$\pi_1(\PPPP^0)=\pi_1(P_t)/<\!\!\tilde R\!\!>_\norm$, where
$\tilde R$ is the set of all the commutators $[c_j,h_k]$, in which
$c_j$ (resp. $h_k$) runs over any set of generators of $\pi_1(P_t)$
(resp. of $\ker [\pi_1(\PPPP^0_U)\rar\pi_1(\PP^{2\dual})]=\pi_1(\PPPP^0_U)$). Using the commutativity
of $\pi_1(P_t)$, as above, we obtain that $[c_j,h_k]=M_\gamma(c_j)c_j^{-1}$,
where $\gamma=f(h_k)\in\pi_1(U)$. Thus for $\PPPP=\PPPP^0$, 
$<\!\!R_1\!\!>_\norm=<\!\!\tilde R\!\!>_\norm=\pi_1(P_t)$,
and we are done.
\end{proof}

\begin{corollary}\label{M-prime}
The partial resolution of singularities $M'$ of $M$ obtained by blowing
up the image of $\Sigma$ is an irreducible symplectic $V$-manifold
whose singularities are $28$ points of analytic type
$ (\CC^4/\{\pm 1\}, 0)$. The natural birational map
$\PPPP^0\dasharrow M'$ is the Mukai flop with center at the
image $\Pi\simeq\PP^2$ of the zero section of $\PPPP^0$, that is, it
blows up $\Pi$ and then blows down the obtained exceptional divisor
$\tilde \Pi\simeq \PP (\Omega^1_\Pi)$ along the second ruling.
The image $\Pi'\simeq\PP^2$ of $\tilde \Pi$ in $M'$ coincides
with the proper transform of $\Sigma_0/\iota_1$.
\end{corollary}

\begin{proof}
To construct $M'$, we may first blow up $\Sigma$ and then quotient by~$\iota_1$.
Let $N=S^{[2]}$, $N_1$ the blowup of $N$ at $\Sigma$, and $N_2$ the blowup of $N_1$ at the  proper transform of $\Sigma_0$. Denote the proper transforms of
$\Sigma$, $\Sigma_0$ in $N_2$ by $\Sigma'$,
$\Sigma'_0$ respectively. The curve of intersection
$\tilde B=\Sigma\cap\Sigma_0$
is a common fixed curve of the pair of regular commuting involutions $\iota_0$,
$\tau$ on $N$, hence it is smooth and
$\Sigma'$, $\Sigma'_0$ intersect transversely along a smooth surface
which is a $\PP^1$-bundle over $\tilde B$.
As the two blowups are done at $\iota_1$-invariant centers, $\iota_1$
lifts to a regular involution, denoted by the same symbol, on $N_2$.
The 3-fold $\Sigma'_0$ and the natural $\PP^1$-bundle
$\Sigma'_0\rar\Sigma_0$ are $\iota_1$-invariant, and $\Fix(\iota_1|_{\Sigma'_0})
=\Sigma'\cap\Sigma'_0$. We deduce that the image $\bar\Sigma'_0$ of
$\Sigma'_0$ in $N_2/\iota_1$ is smooth and is a $\PP^1$-bundle over $\bar\Sigma_0=\Sigma_0/\iota_1$. As we noticed earlier, $\Sigma_0$
is identified with $X$; under this identification, $\iota_1|_{\Sigma'_0}=\iota$,
the Galois involution of $\mu:X\rar\PP^2$. Thus $\bar\Sigma_0\simeq\PP^2$
and $\bar\Sigma'_0\rar\bar\Sigma_0$ is a $\PP^1$-bundle over $\PP^2$.

We have $M'=N_1/\iota_1$, and as the proper transform of $\Sigma_0$
in $N_1$ is isomorphic to $\Sigma_0$, $\bar\Sigma_0\simeq\PP^2$
embeds naturally into $M'$. Denote its image in $M'$ by $\Pi'$.
Then $M''=N_2/\iota_1$ is nothing but the blowup of $M'$ at $\Pi'$,
and we denote by $\tilde\Pi$ the exceptional divisor of this blowup.
The fibers of the blowdown map $\tilde\Pi\rar \Pi'$ represent
one ruling of $\Pi'$, and we are to verify that the map
to $\PPPP^0$ contracts another ruling of $\Pi'$.

Let $\Phi_2:N_2\rar \PPPP^0$ be the composition of $N_2\rar N$
with $\Phi$. The indeterminacy locus of $\Phi$ consists of those
$\xi\in S^{[2]}$ for which $\xi$ is vertical (that is, contained
in a fiber of $\mu\rho$). We omit a fastidious calculation in
local coordinates on $N_2$ which shows that the indeterminacy
is resolved on $N_2$, so that $\Phi_2$ is regular. We can represent
a point $\hat\xi\in N_2$ as a pair $(\xi,\ell_\xi)$, where
$\ell_\xi$ is a line in $\PP^2$ containing $\mu\rho(\xi)$,
and the curve $C_{\hat\xi}=(\mu\rho)^{-1}(\ell_\xi)$ is well defined.
Then $\Phi_2$ contracts to the neutral element $0_{\hat\xi}$
of $\Prym (C_{\hat\xi},\tau)$
all the vertical divisors of $C_{\hat\xi}$. The latter form a curve, isomorphic
to $C_{\hat\xi}/\tau\simeq E_{\hat\xi}:=\mu^{-1}(\ell_\xi)$.
Quotienting further by $\iota_1$, we see that the fiber of the induced
map $\Phi'':M''\rar \PPPP^0$ over $0_{\hat\xi}$ is $E_{\hat\xi}/\iota
\simeq\PP^1$. Thus $\Phi''$ contracts another ruling of $\tilde\Pi$
to the locus $\Pi$ of neutral elements
of the Prymians $P_t$.

As remarked Mukai \cite{Mu-1}, the normal bundle of a plane $\PP^2$
in a symplectic 4-fold is isomorphic to $\Omega^1_{\PP^2}$, so that
the exceptional divisor $\PP(\Omega^1_{\PP^2})$ of the blowup
centered at this $\PP^2$
has exactly two different rulings that can be blown down.
The map $\Phi':M'\rar\PPPP^0$ induced by $\Phi$
blows up $\Pi'$ and contracts the exceptional divisor along
another ruling. It is easily seen, along the lines of the proof
of Lemma \ref{fibers-of-Phi}, that $\Phi'$ is bijective on the
complements to $\Pi,\Pi'$, and this ends the proof.
\end{proof}

We conclude this section by several miscellaneous remarks.

\begin{remark}\label{odd-degree}
{\em Odd-degree Prymians.} It is plausible that all the odd-degree Prymians
$\PPrym^{2k+1,\kappa}(\CCC,\tau)$ contain 3-dimensional rational
subvarieties and thus cannot be symplectic. We will produce such
a subvariety in degree $2k+1=3$ for $\kappa$ defined by $c=C_i'$
as in the paragraph preceding Definition \ref{def-prym}. The case of
degree 3 is particularly handy, because $\PPic^3(|H|)$ is fiberwise
birational to the relative symmetric cube of the linear system $|H|$.
For the fiber $\PPic^3(C)$ corresponding to the reducible
curve $C=C_i'\cup C_i''$, this means that the Abel-Jacobi map
$$AJ:C^{(3)}\dasharrow\PPic^3(C), \ \ \ p_1+p_2+p_3\mapsto
[\OOO_C(p_1+p_2+p_3)]$$ maps birationally all the 4 components
of the symmetric cube $C^{(3)}$ onto the 4 respective components
of $\PPic^3(C)$. We will see that $AJ(C_i^{\prime (2)}
\times C_i'')=\bar J^{2,1}(C)$
is contained entirely in $\PPrym^{3,\kappa}(\CCC,\tau)$.

Let us suppress the subscript $i$ from the notation, so that $c=C',
C=C'\cup C''$. On a typical fiber $C_t$, the involution $\iota$
acts by
\begin{eqnarray*}
\iota :[x_1+x_2+x_3]&\mapsto &[(H+C')\cdot C_t-x_1-x_2-x_3]\\
&=& [(2H-C'')\cdot C_t-x_1-x_2-x_3]\\
&=& [q\cdot C_t - z_1''-z_2''-x_1-x_2-x_3]\\ &=&[y_1+y_2+y_3],
\end{eqnarray*}
where $z_1''+z_2''=C''\cdot C_t$, \ \ $q\in |2H-z_1''-z_2''-x_1-x_2-x_3|$
is a (generically unique) conic passing through the 5 points, and
$y_1+y_2+y_3$ is the residual intersection of this conic with $C_t$.
Now assume $C_t$ \ \ \mbox{$\tau$-invariant}, that is $C_t$ is of the
form $(\mu\rho)^{-1}(\ell_t)$ for a sufficiently general line $\ell_t$.
For generic $x_1,x_2,x_3\in C_t$ (``generic'' here means: which do not
vary in a pencil $g_3^1$), we have:
\begin{eqnarray*}
[x_1+x_2+x_3]\in \PPrym^{3,\kappa}(\CCC,\tau)& \equi&
y_1+y_2+y_3=\tau(x_1+x_2+x_3)\\
& \equi& q\ \mbox{is $\tau$-invariant.}
\end{eqnarray*}
We obtain that the birational transform $\tilde P_t$
of $P_t:=\PPrym^{3,\kappa}(\CCC,\tau)\cap \Pic^3(C_t)$
in $C^{(3)}$ can be described as follows:
$$
\tilde P_t=\{ x_1+x_2+x_3\in C^{(3)}\mid \exists \ q\in\PP^2_{\tau,t}:
\ x_1+x_2+x_3\in q\},
$$
where $\PP^2_{\tau,t}$ denotes the 2-dimensional linear system
of $\tau$-invariant conics 
through the two points   $z_1''+z_2''=C''\cdot C_t$
in the plane spanned by $C_t$.
Now let $\ell_t$ tend to $\ell_0:=\mu\rho(C)$ in the pencil
with a fixed intersection point $p=\ell_0\cap \ell_t$.
Then $z_1''+z_2''$ remains fixed, and the limits of
$\tilde P_t$ contain all the triples $x_1+x_2+x_3$ extracted from
the 6 points $q\cdot C-z_1''-z_2''$, where $q$ runs over the linear
system $\PP^2_{\tau}(z_1'',z_2'')$ of $\tau$-invariant conics
through $z_1'',z_2''$ in the plane $\langle C\rangle$. Varying $p$,
and hence the pair $z_1'',z_2''=\tau(z_1')$, we allow all the
triples $x_1+x_2+x_3$ extracted from the 8-tuples $q\cdot C$,
where $q$ runs over the linear system $\PP^3_\tau$ of all the
$\tau$-invariant conics in $\langle C\rangle$, with the only restriction
that at least one of the two $\tau$-invariant pairs of points
of $q\cap C''$ has empty intersection with $\{x_1,x_2,x_3\}$.
Taking generic points $x_1,x_2\in C'$, $x_3\in C''$, we find  a unique
$\tau$-invariant conic through $x_1,x_2,x_3$, which satisfies the above
restriction, and thus $C^{\prime (2)}
\times C''$ lies in the closure of the family of $\tilde P_t$.
Hence the 3-dimensional rational variety $\bar J^{2,1}(C)$
is contained in $\PPrym^{3,\kappa}(\CCC,\tau)$.
\end{remark}

\begin{remark}
{\em More on the structure of $P_t$.} In Lemma \ref{prym-cc}
and Proposition \ref{prym-strata}, we only enumerated the strata
of the fibers $P_t$; to determine the
topological structure of $P_t$, one should also describe the adjacencies
of these strata. We are going to produce several examples
of such calculation.

In the situation of Lemma \ref{prym-cc}, the open piece $P_0$
of $P_t$ consists of the sheaves 
$$
\FFF=\FFF(0;\lambda_1,\lambda_2,\lambda_3,\lambda_4)=
\OOO_{C_-}\connsum{\lambda_1,\lambda_2,\lambda_3,\lambda_4}\OOO_{C_+}$$
with $\lambda_1\lambda_2=\lambda_3\lambda_4$. Let us fix $\lambda_3,\lambda_4$
and make $\lambda_1\rar 0$; then automatically $\lambda_2\rar \infty$.
The sheaf $\FFF$ can be defined as the subsheaf of $\OOO_{C_-}\oplus\OOO_{C_+}$,
whose stalks at all the points coincide with the stalks of the ambient sheaf,
except at $z_i$, where 
$$
\FFF_{z_i}=\{(f_-,f_+)\in \OOO_{C_-,z_i}\oplus\OOO_{C_+,z_i}\mid \ \ 
f_-(z_i)=\lambda_if_+(z_i)\}.
$$
Thus the stalks of the limiting sheaf $\FFF(0;0,\infty,\lambda_3,\lambda_4)$
coincide with the stalks of $\FFF$ everywhere, except for
the stalks $\mathfrak m_{C_-,z_1}\oplus  \OOO_{C_+,z_1}$ at $z_1$ and
$\OOO_{C_-,z_2}\oplus  \mathfrak m_{C_+,z_2}$ at $z_2$.
Hence
$$
\FFF(0;0,\infty,\lambda_3,\lambda_4)=
\OOO_{C_-}(-z_1)\connsum{\cdot,\cdot,\lambda_3,\lambda_4}\OOO_{C_+}(-z_2),
$$
where the bases of the two sheaves used to define the gluings
at $z_3,z_4$ are the functions $1\in \Gamma(C_\pm,\OOO_{C_\pm})$
considered as rational sections of $\OOO_{C_-}(-z_1)$, $\OOO_{C_+}(-z_2)$.
In the same way, we determine the limit when $\lambda_3,\lambda_4$ are fixed and
$\lambda_1\rar \infty$:
$$
\FFF(0;\infty,0,\lambda_3,\lambda_4)=
\OOO_{C_-}(-z_2)\connsum{\cdot,\cdot,\lambda_3,\lambda_4}\OOO_{C_+}(-z_1).
$$
Changing to the standard bases $e_\pm$ for the sheaves $\OOO_{C_\pm}(-pt)$,
we see that
$$
\FFF(0;0,\infty,\lambda_3,\lambda_4)\simeq \FFF'(0;\frac{z_3-z_2}{z_3-z_1}
\lambda_3,\frac{z_4-z_2}{z_4-z_1}\lambda_4),
$$
$$
\FFF(0;\infty,0,\lambda_3,\lambda_4)\simeq \FFF'(0;\frac{z_3-z_1}{z_3-z_2}
\lambda_3,\frac{z_4-z_1}{z_4-z_2}\lambda_4).
$$
Similarly, we find the limits when $\lambda_1,\lambda_2$ are fixed and
$\lambda_3\rar 0$ or $\infty$. And when $\lambda_1,\lambda_3$ tend simultaneously
to elements of $\{0,\infty\}$, then the limit is the unique sheaf in $P_t$
which is non-invertible at all the 4 points~$z_i$: \ \ 
$\OOO_{C_-}(-2\, pt)\oplus\OOO_{C_+}(-2\, pt)$.
Finally, we conclude:\smallskip 

{\em In the situation of Lemma \ref{prym-cc}, $P_t$ is obtained from
$\PP^1\times\PP^1$ by the following gluings:
\begin{itemize}
\item[--] the horizontal sections $0\times\PP^1$ and $\infty\times \PP^1$ are glued together
according to the rule $(0,\lambda)\sim (\infty, [z_1,z_2;z_3,z_4]^2\lambda)$;
\item[--] the vertical sections $\PP^1\times 0$ and $\PP^1\times \infty$ are glued together
according to the rule $(\lambda,0)\sim ([z_3,z_4;z_1,z_2]^2\lambda,\infty)$;
\item[--] the $4$ ``vertices'' $(0,0),(0,\infty),(\infty,0),(\infty,\infty)$ are
glued together.
\end{itemize}
Here $[z_1,z_2;z_3,z_4]$ stands for the cross ratio of $4$ complex numbers.}\smallskip 

We will also provide the answers for three cases of Proposition \ref{prym-strata},
using the notation introduced in the proof of this proposition.\smallskip 

{\sc Case} (i). {\em The normalization $\tilde P_t$ of $P_t$ is a $\PP^1$-bundle
over the elliptic curve
$E=\Prym (\tilde C,\tau)$ having two distinguished cross-sections $0,\infty$.
Let $0_x$, $\infty_x$ denote the point of the respective cross-section lying
in the fiber over $x\in E$. Then $P_t$ is obtained from $\tilde P_t$
by gluing $0$ to $\infty$ with a translation according to the rule
$0_x\sim \infty_{x+[p''-p']}$.}\smallskip 

{\sc Case} (iii). {\em The normalization $\tilde P_t$ of $P_t$ is a $\PP^1$-bundle
over the elliptic curve
$\tilde C$ having two distinguished cross-sections $0,\infty$,
and $P_t$ is obtained from $\tilde P_t$
by gluing $0$ to $\infty$ with a translation according to the rule
$0_x\sim \infty_{x+[p_1'-p_2'-p_1''+p_2'']}$.}\smallskip 

{\sc Case} (vii). {\em $P_t$ is a locally trivial bundle over the elliptic
curve $\tilde C$ with fiber $\PP^1\bigvee \PP^1$, the bouquet of two copies of
$\PP^1$.}\smallskip 

As concerns the compactified Jacobians of the curves $C_t$, one can find
examples of their calculation in \cite{Cook-2}.
\end{remark}

\begin{remark}
{\em Moduli spaces with involution.} One can pursue our approach
to constructing new symplectic varieties in a generalized setting:
search for pairs $(\MMM,\kappa)$ formed by a moduli space of sheaves on a K3 surface
and a symplectic birational involution. Then one may expect to get new symplectic manifolds
either as a (partial) desingularization of the quotient $\MMM/\kappa$, or
as the fixed locus $\MMM^\kappa$. We can obtain an example of this kind
with $\MMM$ parametrizing non-torsion sheaves by a birational
transformation from the compactified
Jacobian $\PPic^2(|H|)$ of Section \ref{local}. Let $C_i'$ be one of the 56
conics in $S$, $\LLL\in \PPic^2(|H|)$ invertible on its support, and $V=\Ext^1_S
(\LLL\otimes \OOO(-C_i'),\OOO_S)\simeq \CC^2$. Then 
the ext-group classifying the extensions
$$
0\lra V^\dual\otimes \OOO_S\lra\EEE\lra \LLL\otimes \OOO(-C_i')\lra 0
$$
is canonically isomorphic to $\Hom (V,V)$, and
we can define a vector
bundle $\EEE$ as the middle term of this extension
with extension class $\id_V \in \Hom (V,V)$. This provides a birational
isomorphism $\PPic^2(|H|)\dasharrow M_S^{H,ss}(2,H,0)$ in the notation using
the Mukai vector $(2,H,0)=(\rk\EEE,c_1(\EEE),\chi(\EEE)-\rk \EEE)$,
and the (regular) symplectic involution $\kappa$ on $\PPic^2(|H|)$
induces a birational symplectic involution on $M_S^{H,ss}(2,H,0)$.
\end{remark}



\begin{thebibliography}{-----------}

\bibitem[AIK]{AIK} 
Altman, A. B., Iarrobino, A., Kleiman, S. L.:
{\it Irreducibility of the Compactified Jacobian}
In: Real and Complex singularities, Proc. 9th Nordic Summer School
/ NAVF, Oslo, 1976.  Gronignen: Sijthoff \& Noordhoff, 1977, 1--12.

\bibitem[AK]{AK} Altman, A. B., Kleiman, S. L.:
{\it Compactifying the Picard Scheme},
Adv. Math. {\bf 35}, 50--112 (1980).

\bibitem[AKMW]{AKMW} Abramovich, D., Karu, K., Matsuki, K., W\l odarczyk, J.: {\em Torification and factorization of birational maps,} J. Am. Math. Soc. 
{\bf 15}, 531–572 (2002).

\bibitem[Ar]{Ar} Artamkin, I. V.: {\em 
On the deformation of sheaves}, Izv. Akad. Nauk SSSR Ser. Mat.  
{\bf 52}, 660--665, 672  (1988);  
translation in  Math. USSR-Izv.  {\bf 32}, 663--668  (1989).

\bibitem[Art]{Art} Artin, M.: {\em Algebraization of formal moduli, 
I,}  Global Analysis (Papers in Honor of K. Kodaira)  pp. 21--71, 
Univ. Tokyo Press, Tokyo,  1969.

\bibitem[B]{B} Barth, W. {\em Abelian surfaces 
with $(1,2)$-polarization}, In:  Algebraic geometry, Sendai, 1985,  
41--84, Adv. Stud. Pure Math., {\bf 10}, North-Holland, Amsterdam, 1987. 

\bibitem[Beau-1]{Beau-1} Beauville, A.:
{\it Some remarks on K\"ahler manifolds with $c_{1}=0$},
In: Classification of algebraic and analytic manifolds (Katata, 1982), 1--26,
Progr. Math., 39,
Birkhauser:Boston, 1983.

\bibitem[Beau-2]{Beau-2} Beauville, A.: {\it 
Vari\'et\'es k\"ahl\'eriennes dont la
premi\`ere classe de Chern est nulle}, J. Differential Geometry {\bf 18},
755-782 (1983).



\bibitem[Beau-3]{Beau-3} Beauville, A.:
{\it Syst\`emes hamiltoniens compl\`etement int\'egrables associ\'es aux surfaces $K3$},
In: Problems in the theory of surfaces and their classification (Cortona, 1988), 25--31,
SE: Sympos. Math., XXXII,
Academic Press: London, 1991.


\bibitem[BH]{BH}   Brun, J., Hirschowitz, A.: {\em 
Vari\'et\'e des droites sauteuses du fibr\'e instanton g\'en\'eral},
With an appendix by J. Bingener,  
Compositio Math.  {\bf 53}, 325--336  (1984). 

\bibitem[BM]{BM}  Bloch, S.; Murre, J. P.: {\em 
On the Chow group of certain types of Fano threefolds},
Compositio Math. {\bf 39},  47--105 (1979).

\bibitem[CK-1]{CK-1} Choy, J., Kiem, Y.-H.: {\em
Nonexistence of a crepant resolution of some moduli
spaces of sheaves on a K3 surface}, math.AG/0407100.

\bibitem[CK-2]{CK-2} Choy, J., Kiem, Y.-H.: {\em
On the existence of a crepant resolution of some moduli 
spaces of sheaves on an abelian surface},  Math. Z. {\bf 252}, 557--575  (2006).

\bibitem[Cook-1]{Cook-1} Cook, P. R.: {\em
Compactified Jacobians and curves with simple singularities},
Algebraic geometry (Catania, 1993/Barcelona, 1994), 37--47,
Lecture Notes in Pure and Appl. Math., 200,
Dekker, New York, 1998. 

\bibitem[Cook-2]{Cook-2} Cook, P. R.: {\em
Local and global aspects of the module theory of singular curves},
PhD thesis, Univ. of Liverpool, 1993.

\bibitem[D]{D} Debarre, O.: {\em
On the Euler characteristic of generalized Kummer varieties},  
Amer. J. Math. {\bf  121},   577--586  (1999).

\bibitem[Dr]{Dr} Dr\'ezet, J.-M.: {\em
Luna's slice theorem and applications}, In: Algebraic 
group actions and quotients, 39--89, Hindawi Publ. Corp., Cairo, 2004. 



\bibitem[F]{F}  Fujiki, A.: {\em 
On primitively symplectic compact K\"ahler $V$-manifolds of dimension four},
Classification of algebraic and analytic manifolds (Katata, 1982), 
71--250, Progr. Math., 39, Birkhäuser Boston, Boston, MA, 1983.




\bibitem[Gu]{Gu} Gulbrandsen, M. G.: 
{\em Lagrangian fibrations on generalized Kummer varieties}, math.AG/0510145.

\bibitem[Ha]{Ha} Hartshorne, R.:
{\em Residues and Duality,} Lecture Notes in Math. {\bf 20}, Springer-Verlag,
New York, 1966.


\bibitem[HasTsch-1]{HTsch-1} Hassett, B., Tschinkel, Yu.: {\em Abelian 
fibrations and rational points on symmetric products},  
Internat. J. Math. {\bf 11}, 1163--1176 (2000).

\bibitem[HasTsch-2]{HTsch-2} Hassett, B., Tschinkel, Yu.:
{\em Rational curves on holomorphic symplectic fourfolds},  
Geom. Funct. Anal. {\bf 11}, 1201--1228 (2001). 




\bibitem[HL]{HL} Huybrechts, D., Lehn, M.:
{\em The Geometry of Moduli Spaces of Sheaves}, Aspects of Math.,
Vol. E 31, Friedr. Vieweg \&\ Sohn,
Braunschweig (1997). 


\bibitem[Hu-0]{Hu-0} Huybrechts, D.:
{\em  
Birational symplectic manifolds and their deformations},
J. Differ. Geom. {\bf 45}, 488-513 (1997).

\bibitem[Hu-1]{Hu-1} Huybrechts, D.:
{\em   Compact hyper-K\"ahler manifolds: basic results}, 
Invent. Math. {\bf 135},
63--113 (1999); Erratum: {\bf 152},  209--212  (2003).

\bibitem[Hu-2]{Hu-2} Huybrechts, D.:
{\em  The K\"ahler cone of a compact hyperk\"ahler manifold,} 
Math. Ann.  {\bf 326}, 499--513  (2003).

\bibitem[HvM]{HvM} Horozov, E., van Moerbeke, P.: {\em
The full geometry of Kowalewski's top and $(1,2)$-abelian surfaces},
Comm. Pure Appl. Math.{\bf 42}, 357--407 (1989).


\bibitem[IR]{IR}  Iliev, A.,  Ranestad, K.: {\em
The abelian fibration on the Hilbert cube of a K3 surface of genus 9},
e-print math.AG/0507016.




\bibitem[KL]{KL} Kaledin, D., Lehn, M.: {\em
Local structure of hyperkaehler singularities in O'Grady's examples},
math.AG/0405575.

\bibitem[KLS]{KLS} Kaledin, D., Lehn, M., Sorger, C.: {\em
Singular symplectic moduli spaces},  Invent. Math. {\bf 164},
591--614  (2006).

\bibitem[Kir]{Kir} Kirwan, F. C.: {\em
Partial desingularisations of quotients of nonsingular 
varieties and their Betti numbers},
Ann. of Math. (2) {\bf 122}, 41--85 (1985).




\bibitem[Lau]{Lau} Laudal, O. A.: {\em Matric Massey products
and formal moduli, I,}  Algebra, algebraic topology and their 
interactions (Stockholm, 1983),  218--240, 
Lecture Notes in Math., 1183, Springer, Berlin, 1986.

\bibitem[Lei]{Lei} Leibman, A.: {\em Fiber bundles with degenerations and their applications to computing fundamental groups},  Geom. Dedicata  {\bf 48}, 
93--126 (1993).

\bibitem[LP]{LP} Le Potier, J.:
{\em Fibr\'e d\'eterminant et courbes de saut sur
les surfaces alg\'ebriques}, 
Complex projective geometry (Trieste, 1989/Bergen, 1989),  213--240, 
Cambridge Univ. Press, Cambridge, 1992.

\bibitem[LS]{LS} Lehn, M., Sorger, C.: {\em
La singularit\'e de O'Grady}, J. Algebraic Geom. {\bf 15}, 753-770 (2006).

\bibitem[Lu]{Lu} 
Luna, D.:
{\em Slices \'etales}, 
Sur les groupes alg\'ebriques, pp. 81--105. Bull. Soc. Math. France,
Paris, Memoire 33, Soc. Math. France, Paris, 1973. 





\bibitem[Mar-1]{Mar-1} Markushevich, D.: {\em
Some algebro-geometric integrable systems versus classical ones},  
The Kowalevski property (Leeds, 2000),  197--218, CRM Proc. Lecture Notes, 
32, Amer. Math. Soc., Providence, RI, 2002.


\bibitem[Mar-2]{Mar-2} Markushevich, D.: {\em
Rational Lagrangian fibrations on punctual Hilbert schemes of K3 surfaces},
Manuscripta Math. {\bf 120}, 131--150 (2006).

\bibitem[MarT]{MarT} Markushevich, D.,  Tikhomirov, A. S.:
{\it Symplectic structure on a moduli space of sheaves on a cubic fourfold,} 
Izv. Math. {\bf 67}, 121-144 (2003).


\bibitem[Maru]{Maru}  Maruyama, M.: {\em 
Vector bundles on ${\PP}\sp 2$ and torsion sheaves on the dual plane},
In:
Vector bundles on algebraic varieties, Pap. Colloq., Bombay 1984, Stud. 
Math., Tata Inst. Fundam. Res. {\bf 11}, 275-339 (1987). 




\bibitem[Mo]{Mo} Mo{\u\i}{\v{s}}ezon, B. G.:
{\em On {$n$}-dimensional compact complex manifolds having {$n$}
              algebraically independent meromorphic functions, {I},}
	      Izv. Akad. Nauk SSSR Ser. Mat. {\bf 30}, 133--174 (1966).




\bibitem[Mor]{Mor} Morrison, D.: {\em
The geometry of K3 surfaces,} Lectures delivered at the Scuola 
Matematica Interuniversitaria, Cortona, 1988.





\bibitem[Mu-1]{Mu-1}  Mukai, S.: {\em 
Symplectic structure of the moduli space of sheaves 
on an abelian or $K3$ surface},  Invent. Math. {\bf  77},  101--116  (1984).

\bibitem[Mu-2]{Mu-2}  Mukai, S.: {\em On the moduli space of 
bundles on $K3$ surfaces. I}, Vector bundles on algebraic varieties 
(Bombay, 1984), 341--413, Tata Inst.
Fund. Res. Stud. Math., 11, Bombay, 1987. 







 
 

\bibitem[O'G-1]{O'G-1} O'Grady, K. G.:
{\em The weight-two Hodge structure of moduli spaces of sheaves
 on a $K3$ surface},  J. Algebraic Geom. {\bf  6}, 599--644  (1997).

\bibitem[O'G-2]{O'G-2} O'Grady, K. G.:
{\em Desingularized moduli spaces of sheaves on a $K3$},
J. Reine Angew. Math. {\bf 512}, 49--117 (1999).

\bibitem[O'G-3]{O'G-3} O'Grady, K. G.:
{\em A new six-dimensional irreducible symplectic variety},  
J. Algebraic Geom. {\bf  12},  435--505 (2003).


\bibitem[O'G-4]{O'G-4} O'Grady, K. G.:
{\em Irreducible symplectic 4-folds numerically equivalent to $\Hilb^2(K3)$},
math.AG/0504434.

\bibitem[O'G-5]{O'G-5} O'Grady, K. G.:
{\em Irreducible symplectic 4-folds and Eisenbud-Popescu-Walter sextics},
Duke Math. J. {\bf 134}, 99--137  (2006).

\bibitem[P]{P}  Pantazis, S.: {\em Prym varieties and the geodesic 
flow on ${\rm SO}(n)$},  Math. Ann. {\bf  273}, 297--315  (1986).


\bibitem[Rim]{Rim} Rim, D. S.: {\em
Formal Deformation Theory}, SGA 7, Expos\'e VI,
Lect. Notes Math., 288, Springer-Verlag, Berlin-New York, 1972.


\bibitem[S-1]{S-1} Sawon, J.: {\em 
Abelian fibred holomorphic symplectic manifolds},
Turk. J. Math. {\bf 27}, 197-230 (2003). 

\bibitem[S-2]{S-2} Sawon, J.: {\em 
Lagrangian fibrations on Hilbert schemes of points on K3 surfaces},
e-print math.AG/0509224.



\bibitem[SD]{SD} Saint-Donat, B.: {\em Projective models of $K3$ surfaces}, 
Amer. J. Math. {\bf 96}, 602--639 (1974).


\bibitem[Sim]{Sim} Simpson, C. T.: {\em Moduli of representations of the
fundamental group of a smooth projective variety I,} Publ. Math.
I.H.E.S. {\bf 79}, 47--129 (1994).



\bibitem[Y]{Y} Yoshioka, K.: {\em Moduli spaces of stable sheaves on abelian surfaces},  Math. Ann. {\bf 321}, 817--884  (2001). 

\end{thebibliography}
\end{document}